\documentclass [12pt,oneside] {article}
\usepackage{latexsym,amsfonts,amssymb,amscd,amsmath,amsthm,mathrsfs,stmaryrd} 
\usepackage{tikz}

\theoremstyle{plain}
\newtheorem {Thm} {Theorem}[section]
\newtheorem* {Thm*} {Theorem}
\newtheorem* {Prop*} {Proposition}
\newtheorem {Lem}[Thm] {Lemma}
\newtheorem {Prop}[Thm] {Proposition}
\newtheorem {Cor}[Thm] {Corollary}

\theoremstyle {definition}
\newtheorem {Def}[Thm] {Definition}

\newtheorem {Rem}[Thm] {Remark}
\newtheorem {Exa}[Thm] {Example}
\newenvironment{Pf}[1]{{\noindent\sc Proof #1:}}{\qed\\}
\newcommand {\specialmap} [4] {\text {$ #1\negmedspace : #2 #3 #4 $}}
\newcommand {\map} [3] {\specialmap {#1} {#2}{\to} {#3}}
\newcommand {\longmap} [3] {\specialmap {#1} {#2}{\longrightarrow} {#3}}
\newcommand {\isomap} [3] {\specialmap {#1} {#2}{\overset {\cong{\phantom{.}}} 
            {\longrightarrow}} {#3}}

\newcommand {\Hom} {\operatorname {Hom}}

\newcommand {\at}[1] {\arrowvert_{#1}}

\renewcommand {\(} {\left(}
\renewcommand {\)} {\right)} 
\newcommand {\spec} {\operatorname{spec}}
\renewcommand {\geq} {\geqslant}

\newcommand {\sub} {\subseteq}

\newcommand {\CC} {\mathbb C}
\newcommand {\on}[1] {\operatorname{#1}}

\newcommand {\mK} {{\mathcal{K}}}

\newcommand {\pt} {\on{pt}}

\renewcommand {\SS} {\mathbb{S}}

\newcommand {\tensor}{\otimes}

\newcommand {\ZZ} {\mathbb Z}

\newcommand {\AT} {{A_T}}

\newcommand {\Ell} {\mE ll}

\newcommand {\inv} {^{-1}}

\newcommand {\bL} {\mathbb L}

\newcommand {\mC} {\mathcal C}
\newcommand {\mE} {\mathcal E}
\newcommand {\mF} {\mathcal F}
\newcommand {\mH} {\mathcal H}
\newcommand {\mL} {\mathcal L}
\newcommand {\mM} {\mathcal M}
\newcommand {\mO} {\mathcal O}
\newcommand {\mR} {\mathcal R}

\newcommand {\ol}[1]{\overline{#1}}

\newcommand {\RR} {{\mathbb R}}

\usepackage [all]{xy}
\SelectTips{cm}{12}

\begin{document}

\title{The elliptic Weyl character formula}
\author{Nora Ganter\thanks{This work was made possible by a Centenary
  Fellowship from the Faculty of Science at the University of Melbourne
  and by an Australian Research Fellowship.}\\ The University of Melbourne}
\date {\today}

\maketitle
\begin{abstract}
  We calculate equivariant elliptic cohomology of the partial flag
  variety $G/H$, where $H\sub G$ are compact connected Lie groups of
  equal rank. We identify the
  $RO(G)$-graded coefficients $\Ell_G^*$ 
  as powers of Looijenga's line bundle
  and prove that transfer along the map
  $$
    \longmap\pi{G/H}\pt
  $$
  is calculated by the Weyl-Kac character formula. Treating ordinary
  cohomology, $K$-theory and elliptic cohomology in parallel, this
  paper organizes the
  theoretical framework for the elliptic Schubert calculus of
  \cite{Ganter:Ram}. 
\end{abstract}
\section{Introduction}
The topological aspects of representation theory are captured by the
generalized cohomology theory known as equivariant
  $K$-theory. Applied to a point, the $K$-group
\begin{eqnarray*}
  K_G(\pt) & = & R(G)
\end{eqnarray*}
is the representation ring of the structure group $G$. Applied to
other spaces, it yields rings, which are related to the representation
rings. For instance, 
the $K_G$-theoretic 
transfer along the map $\map\pi{G/H}\pt$
gives the induction map
$$\longmap{\on{ind}}{R(H)}R(G).$$
This point of view is taken in \cite{Atiyah:Bott:Weyl_formula}, where
Atiyah and Bott obtain
Weyl's famous formula for the character of $\on{ind}([\varrho])$ 
as an application of their fixed point
formula for the $T$-equivariant transfer $\pi_!$.
Here $T$ is a maximal torus, sitting inside $H$.

Schubert calculus, originally concerned with the 
cohomology of the partial flag varieties,
has long been extended to include the analogous
$K$-theory picture. An essential ingredient in both theories are
pull-backs and transfers (push-forwards) along maps between partial
flag varieties. 

In \cite{Bressler:Evens:Schubert_Calculus}, Bressler and Evens
formulate Schubert calculus in broad generality, replacing cohomology
and $K$-theory with any generalized multiplicative cohomology theory
possessing  
the relevant transfers. The universal example of
such a theory is complex cobordism, and cobordism-theoretic Schubert
calculus is now becoming a discipline of its own (see
\cite{Bressler:Evens:Cobordism}, \cite{Hornbostel:Kiritchenko} and \cite{Calmes:Petrov:Zainoulline}).

We are interested in  
equivariant elliptic cohomology, $\Ell_G$.
It has long been conjectured that 
$\Ell_G$ plays the same role for the representation theory of the loop
group ${\mathcal L G}$ that $K_G$ plays for the representation theory of
$G$. This idea can already be found in Grojnowski's article on the
definition of $\Ell_G$ \cite[p.2 and 3.3]{Grojnowski}, it
took shape in Ando's work on Euler classes \cite{Ando:sigma} and was later picked
up by Lurie \cite{Lurie}.
We will see in Section \ref{sec:coefficients} that 
\begin{eqnarray*}
  \Gamma\Ell^*_G(\pt) & \cong & \widetilde{Th}_*^W
\end{eqnarray*}
is Looijenga's ring of theta functions \cite{Looijenga}. This is where the loop group
characters take their values.\footnote{More precisely, we are
  considering characters of positive energy
  representations of the central extension $\widehat{\mL G}$.} 

The paper at hand is the first in a joint program with Arun Ram,
studying Schubert calculus in elliptic cohomology. For this, the ring
$\widetilde{Th}_*$ will play the same
role as $R(T)$ for the $K$-theoretic Schubert calculus or as the 
symmetric algebra $S(\mathfrak t_\CC^*)$ in
cohomology. 
Our work ties in with the Bressler-Evens program, but
is not a special case of 
\cite{Bressler:Evens:Schubert_Calculus}: Bressler and Evens work Borel
equivariantly, while Grojnowski's $\Ell_T$ is a
genuinely equivariant theory, 
taking values in sheaves over a scheme $\mM_T$. More importantly,
$\Ell_T$ does not possess Thom 
isomorphisms for complex vector bundles, so that the theory of
transfer maps acquires a twist by a line bundle, called the {\em Thom
  sheaf}.

In our main application, the 
Thom sheaf for $\pi_!$ turns out to be $\mL_{Lo}^g$, the
Looijenga line bundle raised to the dual Coxeter number. This accounts
for the shift of level by $g$ occuring in the Weyl-Kac formula.

The paper is organized as follows: treating cohomology, $K$-theory and
elliptic cohomology simultaneously, we view all three theories as
sheaf valued, revisiting, and to some extent reorganizing, the
circle of ideas in 
\cite{Grojnowski}, \cite{Rosu:K-theory}, \cite{Rosu:elliptic},
\cite{Ando:sigma},
\cite{Ginzburg:Kapranov:Vasserot},
\cite{Lurie}, and \cite{Gepner}.

After recalling the definitions and the general setup (Sections
\ref{sec:three_theories} and \ref{sec:chern_construction}), we review
a powerful calculational tool: this is the theory of moment graphs
(Section \ref{sec:moment_graphs}).
We show how to deduce the
isomorphism
\begin{eqnarray*}
  K_T(G/H) &\cong & R(T)\tensor_{R(G)}R(H)
\end{eqnarray*} 
for torsion free $\pi_1G$ (see \cite{McLeod})
directly from the moment graph of $G/H$. Our proof does not
use Pittie's theorem that $R(H)$ is free over $R(G)$, nor does it
involve any explicit calculations with basis elements.
Our argument is identical for $K$-theory and cohomology and also
yields a description of $\Ell_T(G/H)$.   

Now we are in a position to prove the axioms 
of \cite{Ginzburg:Kapranov:Vasserot} needed in our
applications, and this is done in Section \ref{sec:GKV}. 

Section \ref{sec:Thom} treats the theory of Thom sheaves. We identify 
The Thom sheaf $\bL^\xi$ of any complex vector bundle $\xi$ 
with the pull-back of a universal example:
\begin{eqnarray*}
  \bL^\xi&\cong & c_\xi^*\medspace\bL^{univ}.
\end{eqnarray*}
Here $c_\xi$ is the Ginzburg-Kapranov-Vasserot characteristic class
of $\xi$.

Finally, in Section \ref{sec:character_formulas}, we arrive at the promised formula for
the transfer $\pi_!$. In cohomology, this is a formula by Akyildiz and
Carrell \cite{Akyildiz:Carrell}, in $K$-theory, it is the Weyl
formula, and in elliptic cohomology it is the Weyl-Kac formula. 

The combinatorial aspects of the theory, as well as concrete examples
will be addressed in \cite{Ganter:Ram}.
\subsection{Acknowledgments} 
I would like to thank 
Matthew Ando and Ioanid Ro\c{s}u for teaching me much of what I know about
elliptic cohomology.
Many thanks go to Craig Westerland and Alex Ghitza for countless
conversations.
This paper is my account of work that is, to a large extent, joint
with Arun Ram. It is a pleasure to thank him for being such an
inspiring collaborator. 
Finally, I would like to take this opportunity to express my deep gratitude
to Arun, my department and the Faculty of Science for making it
possible to continue this work in spite of very difficult circumstances.
\tableofcontents
\section{The three sheaf valued theories}
\label{sec:three_theories}
Let $G$ be a compact Lie group, and let $X$ be a finite
$G$-CW-complex. We write
$$
  H_G(X) := \sum_{n\in\ZZ}H^{2n}\(EG\times_GX;\CC\)
$$
for (even) Borel equivariant cohomology with complex
coefficients, and
$$
  K_G(X) := \(\on{Vect_G^\CC}(X)\)^{\on{gp}}\tensor_\ZZ\CC
$$
for equivariant $K$-theory (as in \cite{Segal:Equivariant_K-theory}) with complex
coefficients. We will also consider the relative and reduced versions
of these theories. 
These are contravariant functors in $G$ and in $X$ (or pairs $(X,A)$
or $(X,x_0)$).  
For abelian $T$ it follows that 
the coefficient rings 
$$
  H_T:= H_T(\pt)\quad \text{and}\quad K_T:= K_T(\pt)
$$
form Hopf algebras, with the comultiplication given by
multiplication in $T$. For the circle group $U(1)$, we have
\begin{eqnarray*}
  H_{U(1)} &=& \CC[x]\quad\text{and}\\
  K_{U(1)} &=& \CC[z^{\pm1}].
\end{eqnarray*}
These are the Hopf algebras of regular functions on the (affine) group
schemes
\begin{eqnarray*}
\mathbb G_a &=& \mathbb A^1_\CC\quad\hspace{1.05cm}\text{(additive group) and}\\
\mathbb G_m &=& \mathbb A^1_\CC\setminus\{0\}\quad\text{(multiplicative group).}
\end{eqnarray*}
Here
\begin{eqnarray*}
  x &=& c_1(\CC_1)_{U(1)}  
\end{eqnarray*}
is the first Borel equivariant Chern class of the defining
representation $\CC_1$ of $U(1)$. It generates the ideal
\begin{eqnarray*}
    x\CC[x]  & = & I(0)
\end{eqnarray*}
of regular functions on $\mathbb A_\CC^1$ that vanish at $0$. The
$K$-theory class $z$ is the character of $\CC_1$. The
$K_{U(1)}$-theoretic first Chern class of $\CC_1$ equals $1-z$,
generating the ideal $I(1)$ of regular functions on $\mathbb G_m$
vanishing at $1$. 

More generally, let $T$ be a compact abelian Lie group with
Lie algebra $\mathfrak t$, and let 
$$
  \widehat T = \Hom(T,U(1))
$$
be the character lattice of $T$.
Let $T_0\sub T$ be the connected component of $1$, and let
$\Lambda\subset\mathfrak t^*$ be the weight lattice of the torus
$T_0$.
If $T=T_0$ is connected there is an isomorphism
\begin{eqnarray*}
  \Lambda & \stackrel\cong\longrightarrow&\widehat T\\
  \lambda&\longmapsto& e^{2\pi i\lambda}.
\end{eqnarray*}
For $\lambda\in \Lambda$, let $\CC_\lambda$ be the one dimensional
representation of $T_0$ with character $e^{2\pi i\lambda}$. Then the
coefficients $H_T\cong H_{T_0}$ are identified by the Hopf-algebra
isomorphism
\begin{eqnarray*}
  H_{T}&\cong& \Gamma \mathcal O_{\mathfrak t_\CC}\\
  c_1(\CC_{\lambda})_{T_0}&\mapsfrom&\lambda_\CC.
\end{eqnarray*}
Here $\lambda_\CC=\lambda\tensor_\RR\CC$ is
viewed as a regular function 
on the complex algebraic group $\mathfrak t_\CC:=\mathfrak
t\tensor_\RR\CC$. This point of view, going back to Borel, allows us
to interpret $H_T(X,A)$ as the global sections of a coherent sheaf
$\mH_T(X,A)$ on $\mathfrak t_\CC$.

\medskip
In $K$-theory, the $T$-equivariant coefficients are given by the
representation ring
\begin{equation*}
  K_T\medspace =\medspace R(T) \medspace =\medspace  \CC[\widehat T].  
\end{equation*}
For instance,
\begin{eqnarray*}
  K_{T_0}(\pt) & \cong & \CC\{e^\lambda\}_{\lambda\in2\pi i\Lambda},
\end{eqnarray*}
{with}
$e^\lambda e^\mu = e^{\lambda+\mu}$.
So, 
\begin{eqnarray*}
  K_T&\cong& \Gamma\mathcal O_{T_\CC}  
\end{eqnarray*}
is identified with the ring of regular functions of the
complexification $T_\CC$ of $T$. This allows us to view $K_T(X,A)$ as
the global sections of a coherent sheaf $\mK_T(X,A)$ on $T_\CC$.

We have 
\begin{eqnarray*}
  T_\CC &\cong& \Hom\(\widehat T, \CC^\times\), \quad\text{and}\\
  \mathfrak t_\CC& \cong& \Hom(\widehat T, \CC). 
\end{eqnarray*}
Let $E$ be a complex elliptic curve, and let\footnote{$\mM_T$ may be
  interpreted as the moduli scheme of certain principal $T$-bundles
  on $E$, see \cite[(1.4.2)]{Ginzburg:Kapranov:Vasserot}, where $\mM_G$ is
  denoted $\mathcal X_G$.}
\begin{eqnarray*}
  \mM_T &:=& \Hom(\widehat T, E).   
\end{eqnarray*}
Grojnowski's $T$-equivariant elliptic cohomology takes values in
coherent sheaves over $\mM_T$, and 
\begin{eqnarray*}
  \Ell_T(\pt) &= &\mO_{\mM_T}.  
\end{eqnarray*}
We will see that the above theories form the degree zero parts of three
$RO(T)$-graded (sheaf valued) equivariant cohomology theories.
Note that the complex group $\mM_T$ is no longer affine and
the global sections $\Gamma\Ell_T(-)$ do not form a cohomology
theory. This makes the sheaf point of view essential to the theory.
We will now see how the formal properties of $\Ell_T$, as
axiomatically postulated in \cite{Ginzburg:Kapranov:Vasserot},
determine the stalks of the theory. This is the motivation behind
Grojnowski's construction, which we will recall in Section
\ref{sec:construction}. 
\subsection{Homogeneous spaces and representation spheres} 
Let $T'\sub T$ be a closed
subgroup. Then we 
have canonical inclusions 
$\mathfrak t_\CC'\sub\mathfrak t_\CC$ and $T_\CC'\sub T_\CC$
and $\mM_{T'}\sub\mM_T$ and isomorphisms of coherent sheaves
\begin{eqnarray}
\notag
  \mH_T(T/T') & \cong & \mO_{\mathfrak t'_\CC}  \quad\text{(over
    $\mathfrak t_\CC$ )}  \\
\label{eq:homogeneous_spaces}
  \mK_T(T/T') & \cong & \mO_{T'_\CC}  \quad\text{(over
    $T_\CC$)}\\
\notag  \Ell_T(T/T') & \cong & \mO_{\mM_T'}  \quad\text{(over
    $\mM_T$)}.
\end{eqnarray}
The first two of these isomorphisms are classical. We recall the
definition of the third on page
\pageref{page:homogeneous_spaces}. That it is an 
isomorphism will be an immediate consequence of the construction
of $\Ell_T$.

From now on, we let $\AT$ be one of the complex abelian groups
$\mathfrak t_\CC$ or $T_\CC$ or $\mM_T$, and we let $\mF_T$ be the
theory $\mH_T$ or $\mK_T$ or $\Ell_T$ taking values in sheaves over $\AT$.
Often we will write $+$ for the group operation in $A_T$ and $0$ for
its unit, with the understanding that these are to be replaced by
$\cdot $ and $1$ for the multiplicative case $A_T=T_\CC$.

Let $T$ be a torus, $\lambda\Lambda$, and let $\SS^\lambda$
be the representation sphere 
(one point compactification) of $\CC_\lambda$, and write $K_\lambda$
for the kernel of $e^{2\pi i \lambda}$ inside $T$. 
We may identify the
equator of $\SS^\lambda$ with $T/K_\lambda$.
The usual Mayer-Vietoris argument gives the following:
\begin{Cor}
\label{cor:representation_spheres}
  The sheaf
  $\mF_T(\SS^\lambda)$ is identified with the kernel of the map
  \begin{eqnarray*}
    \mO_{A_T} \oplus \mO_{A_T} & \longrightarrow &
    \mO_{A_{K_\lambda}}\\
    (f,g) & \longmapsto & \(f-g\)\at{A_{K_\lambda}}.
  \end{eqnarray*}
\end{Cor}
\subsection{Stalks}\label{sec:stalks}
By a point in $\AT$, we will always mean a maximal point. For
$a\in\AT$, let 
\begin{eqnarray*}
  T(a) &:=& \bigcap_{a\in A_{T'}}T'  
\end{eqnarray*}
be the smallest subgroup of $T$ with $a\in A_{T(a)}$.
Let 
$$
  i_a\negmedspace : X^{T(a)}\hookrightarrow X
$$
be the inclusion of the $T(a)$-fixed points. We will identify the stalk
of $\mF_T$ at $a$ in two steps.

First, we note that 
\begin{equation}
  \label{eq:stalks}
  \isomap{i_a^*}{\mF_T^*(X)_a}{\mF_T^*(X^{T(a)})_a}  
\end{equation}
is an isomorphism of $T$-equivariant cohomology theories. Indeed, it
is enough to check this on orbits $X=T/T'$, where it follows from
\eqref{eq:homogeneous_spaces}.  
Second, consider the quotient map $\map pTT/T(a)$ and use the
isomorphism\footnote{This isomorphism was proved in \cite{Atiyah:Bott:moment} for
  $H$, in \cite{Segal:Equivariant_K-theory} for $K$ and postulated 
for $\Ell$ in \cite[(1.6.3)]{Ginzburg:Kapranov:Vasserot}. Again, it
will follow immediately from the construction of $\Ell_T$.
}
$$
  \mF_T(X^{T(a)}) \cong A_p^*\(\mF_{T/T'}(X^{T(a)})\).
$$
Let 
\begin{eqnarray*}
  \tau_a\negmedspace : \AT & \longrightarrow &\AT \\
   b&\longmapsto & a+b
\end{eqnarray*}
denote translation by $a$. 
Then 
\begin{equation}
  \label{eq:tau_a}
  A_p = A_p \circ \tau_a,  
\end{equation}
and hence
$$
  \mF_T(X^{T(a)})_a \cong   \mF_T(X^{T(a)})_0.
$$
Combining these two steps, we obtain isomorphisms
\begin{eqnarray*}
  \mH_T(X)_a & \cong & H_T(X^{T(a)})\tensor_{H_T}\mathcal O_{\mathfrak
    t_\CC,0} \\
  \mK_T(X)_a & \cong & K_T(X^{T(a)})\tensor_{K_T}\mathcal O_{T_\CC,1}
  \quad\text{and} \\ 
  \Ell_T(X)_a & \cong & \Ell_T(X^{T(a)})_0.
\end{eqnarray*}
Over a sufficiently small 
neighbourhood $U$ of $a$ (see Section \ref{sec:construction} for
details) these isomorphisms extend to an isomorphism of sheaves
\begin{equation}
  \label{eq:neighbourhood}
  \mF_T(X)\at U \quad\!\!\cong\quad\!\! (\tau_{a})_
  *\medspace\(\mF_T(X^{T(a)})\at{U-a}\).   
\end{equation}
\section{The Chern character and the construction of
  $\Ell_T$}
\label{sec:chern_construction}
\subsection{Completion}
In the sheaf-theoretic language, the Atiyah-Segal completion theorem
identifies the formal completion of $\mF_T$ at $0\in A_T$ with the
Borel equivariant version of $\mF$. More precisely, we have the
following theorem.
\begin{Thm}[Completion Theorem]
\label{thm:completion}
  We have an isomorphism of pro-rings
  \begin{eqnarray*}
    \mF_T(X)\widehat{{}_0} &\cong & 
    \mathop{\underleftarrow\lim}_k \medspace\mF(ET^{(k)}\times_T X),
  \end{eqnarray*}
  where $ET^{(k)}$ is the $k$-skeleton of $ET$.
\end{Thm}
In the case of $K$-theory, the right-hand side is $K(ET\times_TX)$,
and Theorem \ref{thm:completion} is 
\cite{Atiyah:Segal:completion}. 
For cohomology, the right-hand side 
is $$\prod_{n\in\ZZ}H^{2n}(ET\times_TX;\CC)$$
(see \cite[p.6]{Rosu:K-theory}).
In Section \ref{sec:Thom},
we will see how Theorem \ref{thm:completion} follows
from the formal properties of $\mF_T$.
\subsection{Ro\c{s}u's Chern character}
\label{sec:Chern} 
Consider the exponential map
$$
  \longmap{\exp}{\mathfrak t_\CC}{T_\CC}.
$$
This is an analytic map of complex groups, it is not
algebraic. 
For an algebraic sheaf $\mF$ on a complex variety, we let $\mF^h$ be
the analytic sheaf associated to $\mF$.
The following theorem is a reformulation of the main result
in \cite{Rosu:K-theory}.
\begin{Thm}[Ro\c{s}u]
  Assume\footnote{This assumption ensures that the restriction maps of
    $\mH_T^h(X)$ and, more importantly, the map
    $\mH_T^h(X)_0\to\mH_T(X)\widehat{{}_0}$ are injective. I do not
    follow Ro\c{s}u's argument for arbitrary $X$ in
    \cite[p.7]{Rosu:K-theory}. This does not affect the main
    applications in \cite{Rosu:K-theory}, since for those Knutsen and
    Ro\c{s}u do require
    $H_T(X)$ to be free over $H_T$.} 
  that $H_T(X)$ is free over $H_T$.
  For a small enough
  analytic open neighbourhood $U$ of $0\in\mathfrak
  t_\CC$ 
  there is an isomorphism of analytic sheaves
  \begin{eqnarray*}
    {ch_T}\negmedspace:{\mK_T^h(X)\at{ \exp(U)}}&\xrightarrow{\phantom{xx}\cong\phantom{xx}}&
    {\exp_*\mH_T^h(X)\at U},    
  \end{eqnarray*}
  uniquely determined by the commuting diagram
  $$
    \xymatrix{
    \mK_T^h(X)_1\ar[0,2]^{ch_T\at1}\ar@{>->}[2,0]&&\mH^h_T(X)_0\ar@{>->}[2,0]\\
    \\
    \mK_T(X)\widehat{{}_1}\ar[0,2]_{ch}&&\mH_T(X)\widehat{{}_0},}
  $$
where the top row is Ro\c{s}u's Chern character at the stalk $1$, and in the bottom
row, $ch$ stands for the (Borel-equivariant) classical Chern character.
\end{Thm}

Consider now the quotient maps
$$\longmap{\exp_E}{\CC}{E = \CC/2\pi i\langle\tau,1\rangle}$$
and 
$$\longmap{y}{\CC^\times}{E \cong \CC^\times/q^\ZZ},$$
where $q= e^{2\pi i\tau}$.
These induce a commuting diagram of complex analytic group homomorphisms
$$
  \xymatrix{
    \mathfrak t_\CC \ar[1,1]_{\exp_{\mM_T}}\ar[0,2]^{\exp_{T_\CC}} && T_\CC
    \ar[1,-1]^{y}\\  
    & \mM_T.
  }
$$
When it exists, Ro\c{s}u's Chern isomorphism $ch_T$ will fit into a commuting diagram
$$
  \xymatrix{
    (\exp_{\mM_T})_*\medspace\mH^h_T(X)\at U && y_*\medspace\!\mK^h_T(X)\at
    {\exp_{T_\CC(U)}}\ar[0,-2]_{y_*(ch_T)}\\
    &\Ell_T^h(X)\at {\exp_{\mM_T}(U)}.\ar[-1,-1]^\phi\ar[-1,1]_\psi
  }
$$
of sheaf isomorphisms over a small neighourhood of $0$ in $\mM_T$.
\subsection{Construction of $\Ell_T$}\
\label{sec:construction}
The construction of $\Ell_T(X)$ was first
outlined in \cite{Grojnowski}. The technical details were filled in in
\cite{Rosu:K-theory}, see also \cite{Ando:sigma} and \cite{Rosu:elliptic}. 
This section is a reminder of Grojnowski's construction.

Note first that the properties stated in Section \ref{sec:stalks}
and Section \ref{sec:Chern} 
determine $\Ell^h_T(X)$ locally: every $a\in\mM_T$ has a small
analytic neighbourhood $U_a$ satisfying
\begin{eqnarray*}
  b\in U_a  &\Longrightarrow &X^{T(b)}\sub X^{T(a)}.  
\end{eqnarray*}
Choose $U_0$ small enough such that the logarithm is 
well-defined over it and
\begin{eqnarray*}
  c\in U_0&\Longrightarrow & X^{T_{\log(c)}} = X^{T(c)}.
\end{eqnarray*}
Further, we assume
$U_a$ to be small enough to satisfy
$$(U_a-a)\sub U_0.$$  
Then we are forced into
\begin{eqnarray*}
  \Ell^h_T(X)\at{U_a} &\cong&
  (\tau_a\circ\exp)_*\medspace\mH^h_T(X^{T(a)})\at{\log(U_a-a)}.
\end{eqnarray*}
We need to understand how these patches are to be glued. 
Given a non-empty intersection $U:=U_a\cap U_b$, we 
make the additional assumptions\footnote{This is possible by
  \cite[2.5]{Rosu:K-theory}.}  
\begin{eqnarray*}
  a-b&\in& U_0,\quad\quad\text{and}\\
  X^{T(b)}&\sub& X^{T(a)}.
\end{eqnarray*}
\begin{Lem}
  Let $i$ denote the inclusion of $X^{T(b)}$ in $X^{T(a)}$. After
  restricting to $\log(U-a)$, the map 
  \begin{eqnarray*}
    i^*\negmedspace :
    {\mH_T(X^{T(a)})}&\stackrel\cong\longrightarrow &
    {\mH_T(X^{T(b)})}
  \end{eqnarray*}
  becomes an isomorphism of (analytic) sheaves.
\end{Lem}
\begin{Pf}{}
  We check the statement on stalks. Let $\gamma\in\log(U-a)$. We claim
  that we have an equality of simultaneous fixed point sets
  $$X^{T(a)}\cap X^{T(\gamma)} = X^{T(b)}\cap X^{T(\gamma)}.$$ 
  The statement then follows from \eqref{eq:stalks} 
  with $\gamma$ in the role of $a$.
  To prove the non-trivial direction of the claim, let
  $c=\exp(\gamma)$. Then $a+c$ is an element of $U_b$, and we obtain 
  $$
    X^{T(a)}\cap X^{T(\gamma)} \quad\! = \quad\!  X^{T(a)}\cap X^{T(c)} 
               \quad\!\sub\quad\! X^{T(a+c)}
               \quad\!\sub\quad\! X^{T(b)}.
  $$
  The inclusion in the middle may be checked on orbits $T/T'\sub X$,
  where it follows immediately from the definition of $T(-)$.
\end{Pf}

Let now $\gamma:=\log(a-b)$.
Similarly to the proof of the lemma, one argues that $T_{\gamma}$
fixes $X^{T(b)}$. As in \eqref{eq:tau_a}, we obtain an isomorphism
\begin{eqnarray*}
  \phi\negmedspace :
    \(\tau_{\gamma}\)_*\medspace\mathcal H_T\(X^{T(b)}\)
   & \cong &   \mathcal H_T\(X^{T(b)}\),
\end{eqnarray*}
and hence of the corresponding analytic sheaves.
Finally, we have
$$
  \tau_b\circ\exp\circ\tau_\gamma = \tau_a\circ\exp.
$$
The desired glueing isomorphism is the composite
$$
  \(\tau_b\exp\)_*(\phi)\circ\(\tau_a\exp\)_*(i^*).
$$
To define the algebraic theory $\Ell_T(-)$, we use 
Serre's GAGA result:
\begin{Thm}[{\cite{Serre:GAGA}}]
  Let $X$ be a projective algebraic variety over $\CC$, let $X^h$ be
  its underlying analytic variety and let $\mC\!oh_{\on{alg}}(X)$ and 
    $\mC\!oh_{\on{an}}\(X^h\)$ be the categories of coherent algebraic
    (resp.\ ananlytic) sheaves over $X$. Then the
  functor
  \begin{eqnarray*}
    \mC\!oh_{\on{alg}}(X) & \longrightarrow &  \mC\!oh_{\on{an}}\(X^h\) \\
                   \mF    &\longmapsto & \mF^h
  \end{eqnarray*}
  is an equivalence of categories.
\end{Thm}
\subsection{Compact connected Lie groups}
Let $G$ be a compact connected Lie group with maximal torus
$T$ and Weyl group $W$. Then $\mF_G$ takes values over the scheme
\begin{eqnarray*}
A_G & = & A_T/W,
\end{eqnarray*}
and 
\begin{eqnarray*}
  \mF_G(X) &=& \mF_T(X)^W  
\end{eqnarray*}
is the sheaf of $W$-invariant sections. In the elliptic case, these
are to be taken as the definition of $\mM_G$ and $\Ell_G$.
\section{Moment graphs}
\label{sec:moment_graphs}
Moment graph theory provides a powerful tool for calculations.
Let $T$ be a compact torus and $X$ a compact $T$-manifold. 
Let $$\longmap i{X^T}X$$ be the
inclusion of the fixed points in $X$. For a subtorus $T'$ of $T$, this
factors through the inclusion
$$\longmap{i_{T'}}{X^T}{X^{T'}}.$$

Recall that
the equivariant 1-skeleton $X_1$ of $X$ is defined as the set of all
points in $X$ whose orbit is at most one-dimensional.
The following theorem was proved for cohomology by Goresky, Kottwitz and
MacPherson \cite{Goresky:Kottwitz:MacPherson}. Knutsen and Ro\c{s}u later
generalized it  
to $K$-theory and elliptic cohomology
\cite{Rosu:K-theory}.
\begin{Thm}[Localization Theorem]\label{thm:localization}
Assume that $H_T(X)$ is free over $H_T$ and that $X_1$ consists of a
finite number of representation spheres 
$\SS^{\lambda}$, meeting only at the fixed points.
Then the map
\begin{eqnarray*}
  i^*\negmedspace : \mF_T(X) & \longrightarrow & \mF_T(X^T)
\end{eqnarray*}
is injective, and its image is equal to 
\begin{equation}
  \label{eq:localization}
  \on{Im}(i^*) = \bigcap_{T'} \on{Im}(i_{T'}^*),  
\end{equation}
where the intersection runs over all subgroups of codimension
$1$ in $T$.
\end{Thm}
The data determining the right-hand side of \eqref{eq:localization}
are recorded in the ``moment graph'' of $X$:  
\begin{Def}
  In the situation of the theorem,
  the {\em moment graph} $\Gamma$ of $X$ has vertices indexed by the
  fixed points of $X$ and an oriented edge with label $\lambda\in\Lambda$ from
  $x_1$ to $x_2$ for each $\SS^\lambda\sub X_1$, containing $x_1$ as $0$ and
  $x_2$ as $\infty$.
\end{Def}
\begin{Cor}[{of Theorem \ref{thm:localization} and Corollary
    \ref{cor:representation_spheres}}]\label{cor:moment} 
  In the situation of the theorem, $\mF_T(X)$ is 
  described by the following equalizer diagram
  $$
    \mF_T(X)\medspace\longrightarrow\medspace
    \bigoplus_v\mO_{A_T} \medspace
     \tikz[minimum height=0ex]
     \path[->, thick]
     node (a)            {}
     node (b) at (1.5em,0) {}
     (a.north)  edge (b.north)
     (a.south)  edge (b.south);
     \medspace
    \bigoplus_{(e,\lambda)} \mO_{A_{K_\lambda}}.
  $$
  Here $K_\lambda=\ker(e^{2\pi i\lambda})$, 
  the first sum is over the vertices, the second sum is over
  the edges of the moment graph, and the two arrows are defined in the
  obvious manner.
\end{Cor}
This formulation of the theory can be found in the paragraph before (1.3) in
\cite[p.27]{Goresky:Kottwitz:MacPherson}. 
%
\begin{Exa}[Partial flag varieties]
\label{exa:flag}
  Let $H\sub G$ be compact connected Lie groups of equal rank. Let
  $T\sub H$ be a maximal torus (of both). Let $W_G$ be the Weyl group of
  $G$. Then the Weyl group of $H$ can be identified with a subgroup
  $W_H\sub W_G$. Fix a set of positive roots $\mR_+$ of
  $G$. 
  For $\alpha\in\mR_+$, let $s_\alpha\in W_G$ be the corresponding reflection. 
  The following description of the 
  moment graph of $G/H$ can be found in 
  \cite[Thm 3.1]{Tymoczko}:
  we have a bijection
  \begin{eqnarray*}
     W_G/W_H & \longrightarrow &(G/H)^T\\
     wW_H  & \longmapsto & wH.
  \end{eqnarray*}
  Writing $[w]$ for the vertex corresponding to the left-coset $wW_H$, 
  we have an edge labeled $\alpha$ 
  from $[w]$ to $[s_\alpha w]$ whenever $\alpha \in\mathcal
  R_+$ is such that $w\inv(\alpha)<0$ is not a root of $H$. 
\end{Exa}
Often the groups $K_\alpha =
\ker(e^{2\pi i\alpha})$, turning up as the stabilizers of one
dimensional orbits in $G/H$, 
have an interpretation as fixed points: 
\begin{Lem}\label{lem:torsion-free}
  Let $G$ be a compact connected Lie group with maximal torus $T$ and
  Weyl group $W$. Let $\alpha$ be a root of $G$. Then the action of
  $s_\alpha$ on $T$ leaves the elements of $K_\alpha$
  fixed. If $\pi_1(G)$ is torsion free then the inclusion
  \begin{eqnarray*}
    K_\alpha & \sub & T^{s_\alpha}
  \end{eqnarray*}
  is an equality.
\end{Lem}
\begin{Pf}{}
  The first claim is \cite[V.(2.9)(iii)]{Broecker:tomDieck}. Recall
  from \cite[V.(7.1)]{Broecker:tomDieck} that 
  \begin{eqnarray*}
    \pi_1(G) &=& \Lambda\!^{\tiny\vee}/\Gamma,
  \end{eqnarray*}
  where $$\Lambda\!^\vee = \ker(\exp)\sub\mathfrak t$$
  and $\Gamma$ is the sublattice
  generated by the coroots. Let $x\in\mathfrak t$ be such that
  $\exp(x)$ is fixed under $s_\alpha$. Then 
  \begin{eqnarray*}
    \alpha(x)\check\alpha & = & x - s_\alpha(x)    
  \end{eqnarray*}
  is an element of $\Lambda\!^\vee$. Since
  $\Lambda\!^{\tiny\vee}/\Gamma$ is torsion free, it follows that
  $\alpha(x)$ is an integer. Hence $e^{2\pi i \alpha(x)} = 1$, and we
  have proved 
  \begin{eqnarray*}
    \exp(x) & \in & K_\alpha.
  \end{eqnarray*}
\end{Pf}
\subsection{The scheme $X_{A_G}$}
\label{sec:X_A}
The sheaf $\mF_G(X)$ is a sheaf of commutative algebras over $A_G$. 
Following \cite[(1.7.4)]{Ginzburg:Kapranov:Vasserot}, we let $X_{A_G}$
be the spectrum of $\mF_G(X)$. This is a scheme over $A_G$. The
assignment $$X\longmapsto X_{A_G}$$ is covariantly functorial in $X$. We
have $\pt_{A_G}=A_G$. 
Writing $\map\pi X\pt$ for the unique map from $X$ to the one point space, the map
$$
  \longmap{\pi_{A_G}}{X_{A_G}}{A_G}
$$
is the structure morphism.
In other words, $X_{A_G}$ is determined by the fact that
\begin{eqnarray*}
  \(\pi_{A_G}\)_*\mathcal O_{X_{A_G}} &\cong& \mF_G(X).  
\end{eqnarray*}
\subsection{Partial flag varieties}
Let $H\sub G$ be compact, connected Lie groups of equal rank, let
$T\sub H$ be a maximal torus, $W_G$ the Weyl group of $T$ in $G$.
\begin{Thm}\label{thm:flag}
  Assume that $\pi_1(G)$ is torsion free. Then
  we have a $W_G$-equivariant epimorphism of schemes over $A_T$
  \begin{eqnarray*}
   \varphi\negmedspace:
   (G/H)_{A_T} &\longrightarrow &
    {A_T\times_{A_G}A_H}, 
  \end{eqnarray*}
  inducing an isomorphism of schemes over $A_G$
  \begin{eqnarray*}
    (G/H)_{A_G} &\cong & {A_H}.
  \end{eqnarray*}
\end{Thm}

In the case of ordinary cohomology, the assumption that $\pi_1(G)$ is
torsion free is not needed.

\medskip
\begin{Pf}{of Theorem \ref{thm:flag}}
  With the notation as in Example \ref{exa:flag}, we write $$F :=
  W_G/W_H$$
  for the set of vertices in the moment graph. 
  Consider the map
  \begin{eqnarray*}
   \ol\varphi\negmedspace : \coprod_{[w]\in F} A_T & \longrightarrow&
   A_T\times_{A_G}A_H\\ 
           ([w],a) & \longmapsto & (a,[w\inv a]).
  \end{eqnarray*}
  Let $W_G$ act on the source of $\ol\varphi$ by
  \begin{eqnarray*}
    v\cdot([w],a) &=& ([vw],va),
  \end{eqnarray*}
  and on the target by its usual action on the first factor. Then
  $\ol\varphi$ is $W_G$-equivariant. 
  By Corollary \ref{cor:moment} and Example \ref{exa:flag}, we have a
  coequalizer diagram 
  \begin{eqnarray*}
    \coprod_{[w],\alpha} A_{K_\alpha} &{     \tikz[minimum height=0ex]
     \path[->, thick]
     node (a)            {}
     node (b) at (1.5em,0) {}
     (a.north)  edge (b.north)
     (a.south)  edge (b.south);} &
    \coprod_{[w]\in F} A_T \medspace\longrightarrow \medspace(G/H)_{A_T}\\ 
    ([w],\alpha,a) & \longmapsto& ([w],a)\\
    ([w],\alpha,a) & \longmapsto& ([s_\alpha w],a).
  \end{eqnarray*}
  where the first coproduct runs over the edges of the moment
  graph. Both maps on the left are $W_G$-equivariant with respect to the action
  \begin{eqnarray*}
    v\cdot([w],\alpha,a) &=& ([vw],v(\alpha),va)
  \end{eqnarray*}
  on their source.
  The universal property of coequalizer yields a $W_G$-equivariant map
  \begin{eqnarray*}
    \varphi\negmedspace : (G/H)_{A_T} &\longrightarrow & A_T\times_{A_G}A_H,
  \end{eqnarray*}
  which is easily seen to be an epimorphism.
  Note that we have identified $\ol\varphi$
  with $A_i$, where $i$ is the inclusion of the $T$-fixed points $F$ in $G/H$. 
  To obtain the promised map of schemes over $A_G$, we
  quotient by the action of $W_G$.
  It remains to prove injectivity of $\varphi/W_G$. 
  We have 
  \begin{eqnarray*}
    \ol\varphi([w],a) = \ol\varphi([v],b) &\iff& a=b \quad\text{ and
    }\quad [w\inv a] = [v\inv a].
  \end{eqnarray*}
  In the source of $\varphi$, we have made the, a priori finer,
  identifications
  \begin{eqnarray}\label{eq:reflections}
    ([w],a) \sim ([s_\alpha w],a) &:\iff& s_\alpha a = a.
  \end{eqnarray}
  Here we have used Lemma \ref{lem:torsion-free}, which is why we need
  the assumption that
  $\pi_1 (G)$ be torsion free.
  In many cases \eqref{eq:reflections} is sufficient to imply injectivity
  of $\varphi$, but we will see an example where this fails.
  Assume now that $[w\inv a]=[v\inv a]$. Then there is an element
  $u\in W_H$ with $$w\inv a=uv\inv a.$$
  In
    $(G/H)_{A_G}$ 
  we have
  \begin{equation*}
    ([w],a)  \sim  ([1],w\inv a) 
             =  ([1],uv\inv a) 
             \sim  ([vu\inv],a)
             =  ([v],a).    
  \end{equation*}
  Hence $\varphi/W_G$ is an isomorphism, as claimed.
\end{Pf}

We now ask when the map $\varphi$ of the Theorem is injective.
\begin{Lem}\label{lem:fixed_points}
  Let $w\in W_G$ and let $T^w_c$ be a connected component of the
  subgroup of $w$-fixed points in $T$.
  Then we can write $w$ as a word 
  \begin{eqnarray*}
    w &=& s_{\alpha_{1}}\cdots s_{\alpha_{l}}
  \end{eqnarray*}
  in (not necessarily simple) reflections such that
  \begin{eqnarray*}
    {T}_c^w & \sub & T^{s_{\alpha_1}}\cap\dots\cap T^{s_{\alpha_l}}.
  \end{eqnarray*}
\end{Lem}
\begin{Pf}{}   
  Choose $t\in T_c^w$ with $T_c^w\sub\ol{\langle t\rangle}$. Let
  $Z_G(t)$ be the centralizer of $t$ in $G$. This is a connected
  closed subgroup of full rank.
  Its Weyl group $W_Z$ may be viewed as a
  reflection subgroup of $W_G$. All elements of $W_Z$ fix $t$ and
  hence $T^w_c$.
  Since $w\in W_Z$, we can write $w$
  as a word in the reflections generating $W_Z$.
\end{Pf}

The following example shows that the $s_{\alpha_j}$ in Lemma
\ref{lem:fixed_points} can not always be
chosen independently of the connected component.
\begin{Exa}
  Let $G=G_2$, and consider the element $w\in W$ acting by $(-)\inv$ on
  $T$. Then $T^w = T[2]$ has four elements. For each non-trivial element
  of $t\in T^w$ there is a different, unique pair of reflections
  $s_{\alpha_t},s_{\beta_t}\in W$ fixing $t$. For each such pair, 
  $w=s_{\alpha_t} s_{\beta_t}$.   
\end{Exa}
\begin{Cor}[{of Lemma {\ref{lem:fixed_points}}}]
  In the situation of the lemma, we have 
  \begin{eqnarray*}
     \mathfrak t_\CC^w & = & \mathfrak t_\CC^{s_{\alpha_1}}\cap\dots\cap \mathfrak
     t_\CC^{s_{\alpha_l}} 
  \end{eqnarray*}
and
  \begin{eqnarray*}
     \(T_\CC^w\)_c & \sub & T_\CC^{s_{\alpha_1}}\cap\dots\cap
     T_\CC^{s_{\alpha_l}}. 
  \end{eqnarray*}
\end{Cor}
\begin{Cor}
  For cohomology and $K$-theory, the map $\varphi$ of Theorem
  \ref{thm:flag} is an 
  isomorphism. If the centralizers of commuting pairs in $G$ are
  connected, then $\varphi$ is also an isomorphism in elliptic cohomology.
\end{Cor}
On global sections, $\varphi$ gives the familiar isomorphisms
\begin{eqnarray*}
  H_T\tensor_{H_G}H_H & \cong &   H_T(G/H), 
\end{eqnarray*}
studied by Borel, Demazure and others, and
\begin{eqnarray*}
  R(T) \tensor_{R(G)} R(H) & \cong & K_T(G/H)
\end{eqnarray*}
\cite{McLeod}.
\begin{Rem}
  The condition that the centralizers of
  commuting pairs be connected, should be compared to
  \cite[3.2]{Grojnowski}.
\end{Rem}
\begin{Exa}
  Let $G=U(n)$. Then all the groups $T^w$ are, in fact, connected, so
  that the word in Lemma \ref{lem:fixed_points} depends only on
  $w$, so that the map $\varphi$ of the theorem
  is an isomorphism also in the elliptic case.
\end{Exa}
\section{The Ginzburg-Kapranov-Vasserot characteristic class} 
\label{sec:GKV}
\subsection{Properties of $X_{A_G}$}
We discuss some basic properties of the scheme $X_{A_G}$ introduced in
Section \ref{sec:X_A}.
In the cases of
cohomology and $K$-theory, these are well-known. In the elliptic case,
they were conjectured 
in \cite{Ginzburg:Kapranov:Vasserot}. We give proofs for the special
cases relevant to us.

\medskip
\noindent
{\bf Change of groups:} Let $\map\phi HG$ be a map of groups,
and let $X$ be a finite $G$-CW-complex. Then we have a
commuting square
 $$
    \xymatrix{
      X_{A_H}\ar[r]\ar[d]_{\pi_{A_H}}  &
      X_{A_G}\ar[d] ^{\pi_{A_G}}\\
      A_H  \ar[r]_{A_\phi} &
      A_G,
    }
  $$
  We write $X_{A_\phi}$ for the top map. The assignment $(-)_{A_\phi}$
  is natural in $X$.

\medskip
\noindent
{\bf Induction Axiom:} 
Let $K\vartriangleleft G$ be a normal
  subgroup, and let $X$ be a $G$-space such that the action of $K$ on
  $X$ is free. Write $\map p XK\backslash X$ and $\map\phi GG/K$ for
  the quotient maps. Then we have a 
  commuting square, natural in the space $X$,
  $$
    \xymatrix{
      X_{A_G}\ar@{=}[r]^{\sim\phantom{xxX}}\ar[d] &
      \(K\backslash X\)_{A_{G/K}} \ar[d] \\
      A_G\ar[r]_{A_\phi} &
      A_{G/K}, 
    }
  $$ 
  where the vertical maps are the respective structure maps, and the
  top map is $(K\backslash X)_{A_\phi}\circ p_{A_G}$.

\medskip
\noindent
{\bf Homogeneous spaces:} 
  \label{page:homogeneous_spaces}
  Let $j\negmedspace : H\hookrightarrow
  G$ be the inclusion of a 
  closed subgroup. Let $\map\iota\pt G/H$ denote the inclusion of the point
  $1H$.
  Then we have an isomorphism $$I^G_H\negmedspace :(G/H)_{A_G}\cong A_H,$$
  fitting into the commuting diagram
  $$
    \xymatrix{
      & (G/H)_{A_H} \ar[1,-1]_{(G/H)_{A_j}} & \\%
      \(G/H\)_{A_G}
      \ar[1,1]_{\pi_{A_G}}\ar@{=}[0,2]^\sim_{I^G_H}&& 
      A_H \ar[-1,-1]_{\iota_{A_H}}
      \ar[1,-1]^{A_j}\\ 
      & A_G .&
    }  
  $$

\medskip
\noindent 
{\bf K\"unneth:}
For a $G$-space $X$ and an $H$-space $Y$, we have a commuting square,
natural in all ingredients,
  $$
    \xymatrix{
      \(X\times Y\)_{A_{G\times H}}\ar@{=}[r]^{\phantom{x}\sim}\ar[d] &
      X_{A_G}\times Y_{A_H} \ar[d] \\
      A_{G\times H}\ar@{=}[r]^{\sim\phantom{X}} &
      A_G\times A_H,
    }
  $$
  whenever the corners are defined.
  In the special case that $G=H=T$ is a compact torus, we have an
  isomorphism over $A_T$
  \begin{eqnarray*}
      \(X\times Y\)_{A_{T}} &\cong &   X_{A_T}\times_{A_T} Y_{A_T}.
  \end{eqnarray*}

Each of these properties can be reformulated in terms of the sheaves
$\mF_G(X)$, where the obvious generalization for pairs can be
formulated (see \cite{Ginzburg:Kapranov:Vasserot}). 
The last property that we need is stated most naturally in terms of
the reduced theory.

\medskip
\noindent 
{\bf Odd coefficients:}
Let $\varrho$ be an odd-dimensional orthogonal representation of
$G$. Then $\widetilde\mF_G$ vanishes on the corresponding
representation sphere:
\begin{eqnarray*}
  \widetilde\mF_G\(\SS^\varrho\) & = & \{0\}.
\end{eqnarray*}
For elliptic cohomology, 
the change of groups
property, the K\"unneth property and the vanishing of the odd
coefficients
follow immediately from the construction of $\Ell_G$
and the corresponding properties of $H_G$. 
\begin{Prop}\label{prop:equivalent_properties}
  Assuming the change of groups and K\"unneth properties,
  the homogeneous spaces property and the induction axiom are equivalent. 
\end{Prop}
\begin{Pf}{}
  It is shown in \cite[(1.7.5)]{Ginzburg:Kapranov:Vasserot} that the
  induction axiom
  implies the homogeneous spaces property. The other
  direction is proved by cellular induction: if $K$ acts
  freely on the orbit $G/H$ then the composite
  $$
    H\hookrightarrow G \rightarrow
    K\backslash G
  $$ 
  is still injective, and we have
  \begin{eqnarray*}
    (K\backslash G/H)_{A_\phi}\circ p_{A_G}\circ I^G_H & = &
    I^{K\backslash G}_H.
  \end{eqnarray*}
Hence $(K\backslash G/H)_{A_\phi}\circ p_{A_G}$ is an isomorphism.
\end{Pf}{}

We saw in \eqref{eq:homogeneous_spaces} that the homogeneous spaces property 
holds if $G$ is a compact torus.
It follows that the induction axiom holds for the inclusion
$H\sub T$ of any closed subgroup of a compact torus. In particular,
  $G_{A_T}  =  \(T\backslash G\)_{A_1}$,
and hence
\begin{eqnarray}
  G_{A_G} & = & \spec(\CC). \label{eq:G_A_G}
\end{eqnarray}
Further, we saw in 
Theorem \ref{thm:flag} that the homogeneous spaces property holds if 
$H\sub G$ are compact and
connected of equal rank and $\pi_1(G)$ is torsion free.
\begin{Prop}
\label{prop:induction_axiom}
  Let $G$ and $K$ be compact connected Lie groups. Then the induction
  axiom holds for the inclusion $$K\vartriangleleft G\times K.$$ 
\end{Prop}
\begin{Pf}{}
  Write $T$ and $T'$ for the maximal tori of $G$ and $K$. Using
  K\"unneth, \eqref{eq:homogeneous_spaces} and \eqref{eq:G_A_G}, we see that the
  homogeneous spaces property holds for any inclusion of the form 
  $$
    j\times 1\negmedspace : H\times \{1\}\longrightarrow T_G\times K.
  $$
  The case $G=T_G$ now follows from the 
  proof of Proposition \ref{prop:equivalent_properties}.
  In the general case, we have
  \begin{eqnarray*}
    X_{A_{G\times K}} 
    &\cong & \(X_{A_{T\times T'}}\)/\(W_G\times W_K\) \\
    &\cong & \(X_{A_{T\times K}}\)/W_G \\
    &\cong & \(\(K\backslash X\)_{A_{T}}\)/W_G\\
    &\cong & \(K\backslash X\)_{A_{G}}.
  \end{eqnarray*}
This completes the proof.
\end{Pf}{}
\subsection{Classifying maps}
Let $X$ be a compact $G$-manifold, and let 
$$
  \longmap\xi PX
$$
be a $G$-equivariant principal $K$-bundle on $X$. We make the
convention that both groups act from the left and that the actions
commute.
 
Write $G\ltimes X$ for the translation groupoid with objects $X$,
arrows
$$
  x \xrightarrow{\phantom{X}g\phantom{X}} gx
$$
and composition given by composition in $G$. Similarly, we
have the translation groupoids $(G\times K)\ltimes P$ and $K\ltimes\pt$.
\begin{Def}
\label{def:classifying_map}
  The {\em classifying map} of $\xi$ is the generalized map of Lie
  groupoids 
  $$
    f_\xi\negmedspace : G\ltimes X\xleftarrow{\phantom{XX}\simeq\phantom{XX}}
    (G\times K)\ltimes P \xrightarrow{\phantom{XXxXX}} K\ltimes \pt.
  $$
\end{Def}
\begin{Def}\label{def:universal_bundle} 
  The {\em universal principal $K$-bundle} is the $K$-equivariant
  principal bundle
  $$
    \longmap{\xi_{\on{univ}}}{K\ltimes K}{K\ltimes\pt},
  $$
  where the two left-actions of $K$ are as follows:
  as an $K$-equivariant space, $K$ carries
  the action of $K$ on itself by left-multiplication. 
  It is a principal $K$-bundle over the one point
  space via the action $(k_1,k_2)\mapsto k_2k_1\inv$.
\end{Def}
The nomenclature is justified by the following Lemma, which follows
directly from the defintions.
\begin{Lem}
  We have an isomorphism
  \begin{eqnarray*}
    f^*_\xi\(\xi_{\on{univ}}\)&\cong&\xi.    
  \end{eqnarray*}
\end{Lem}
More can be said here: the assignement
$$
  \xi\longmapsto f_\xi
$$
wants to be an equivalence from the category of principal $K$-bundles
over $G\ltimes X$ to the category of generalized maps (i.e., zig-zags
like the one in Definition
\ref{def:classifying_map}), and in fact the former has been used to define
the latter, see \cite{Lerman} and \cite{Hilsum:Skandalis}.  
\begin{Cor}
  In the situation of Definition \ref{def:classifying_map}, assume
  that $G$ is trivial. Then the Borel construction functor, applied to
  $f_\xi$,returns the zig-zag 
  $$
    X\xleftarrow{\phantom{X}\simeq\phantom{X}}EK\times_KP\xrightarrow{\phantom{XXX}}
    BK.
  $$
  Choosing a homotopy inverse to the first map, we obtain the more
  familiar classifying map from $X$ to the classifying space of $K$.
\end{Cor}
\begin{Pf}{}
  This follows, since Borel construction commutes with pull-backs.
\end{Pf}{}
\begin{Exa}[Representations]\label{exa:representations}
  Let $\map\varrho GU(n)$ be a complex representation of $G$, and
  let $\xi_n$ be the universal principal $U(n)$-bundle
  as in Definition \ref{def:universal_bundle}. Consider the action of
  $G$ on $U(n)$ by left multiplication with $\varrho(g)$. This makes
  $\xi_n$ into a $G$-equivariant principal $U(n)$-bundle over the one
  point space. The equivalence 
  \begin{eqnarray*}
    G\ltimes\pt & \stackrel\simeq\longleftarrow & (G\times U(n)) \ltimes U(n)
  \end{eqnarray*}
  of Definition \ref{def:classifying_map} has the quasi-inverse
  \begin{eqnarray*}
    g & \longmapsto & (g,\varrho(g)) \in \on{Stab}(1). 
  \end{eqnarray*}
  Hence the classifying map $f_{\xi_n}$ is equivalent to 
  $$
    \longmap\varrho{G\ltimes \pt}{U(n)\ltimes\pt}.
  $$
\end{Exa}
\begin{Exa}[The splitting principle]
Assume that $K$ is a compact connected Lie
group and $i\negmedspace : T\hookrightarrow K$  the inclusion of its maximal 
torus. Then $\xi$ may be factored as the composite
$$
  \xi\negmedspace : P\xrightarrow{\phantom{X}\zeta\phantom{X}} T\backslash P
  \xrightarrow{\phantom{X}q\phantom{X}} X,
$$
where $\zeta$ is the quotient map by the $T$-action. So, $\zeta$ is a
principal $T$-bundle, while the fiber of $q$ is the flag variety
$T\backslash K$. Over the total space $T\backslash P$ of this flag
bundle, the structure group of $\xi$ can be  
reduced to $T$. 
Let
$$
  \longmap{\zeta[K]}{K\times_TP}{T\backslash P}
$$
be the principal $K$-bundle obtained from $\zeta$ by associating the
fiber $K$.
Then
\begin{eqnarray*}
  q^*(\xi) &\cong & \zeta[K].
\end{eqnarray*}
This fact is known as the {\em splitting principle}.

In terms of classifying maps, the splitting principle amounts to
the commutativity of the diagram 
$$
  \xymatrix{
  f_\zeta\negmedspace: G\ltimes (T\backslash P) \ar[d]_q  &
  (G\times T) \ltimes P \ar[l]_{\simeq}\ar[r]\ar[d] &
  T\ltimes\pt \ar[d]^i\\
  f_\xi\negmedspace :G\ltimes X  &
  (G\times K) \ltimes P \ar[l]_{\simeq\phantom{x}}\ar[r] &
  K\ltimes\pt.
  }
$$
\end{Exa}
\begin{Def}[Characteristic class] 
\label{def:characteristic_class}
Let $\map\xi PX$ be a $G$-equivariant principal $K$-bundle with
classifying map $f_\xi$. Then Proposition \ref{prop:induction_axiom} yields
a map of schemes
$$
  \xymatrix{
  c_\xi\negmedspace: X_{A_G}  &
  P_{G\times K} \ar@{=}[l]\ar[r] &
  A_K.
  }
$$  
We will refer to $c_\xi$ as the {\em Ginzburg-Kapranov-Vasserot
  characteristic class} of $\xi$.
\end{Def}
\section{Thom sheaves}
\label{sec:Thom}
Let $\map\xi VX$ be a $G$-equivariant complex vector bundle. We will
write $X^\xi$ for the Thom space of $\xi$ and $z\negmedspace:
X_+\hookrightarrow X^\xi$ for the zero section. Applying the reduced
theory, we obtain a locally free rank one\footnote{For cohomology and
  $K$-theory this is a classical result.
  In the elliptic
  case it is an immediate consequence of \cite[2.6]{Grojnowski} and the
$W$-equivariance of the cohomology Thom isomorphism.} module sheaf $\widetilde
\mF_G\(X^\xi\)$ over $\mF_G(X)$.
\begin{Def}
  The {\em Thom sheaf}
  of $\xi$ is the line bundle  $\bL_G^{\xi}$ over
  $X_{A_G}$ characterized (up to isomorphism) by
  $$
    \pi_{A_G*}\(\bL_G^{\xi}\)\inv \cong \widetilde\mF_G\(X^\xi\).
  $$
  Note that our convention differs from that in
  \cite[2.1]{Ginzburg:Kapranov:Vasserot}, where the inverse of
  $\bL^\xi_G$ is refered to as the Thom sheaf. 
  The {\em Euler map} is the map 
  $$
    \longmap{\eta_G^\xi}{\mathcal O_{X_{A_G}}}{\bL_G^\xi}
  $$
  induced by the zero section $\map zXX^\xi$ (compare
  \cite[(2.6)]{Ginzburg:Kapranov:Vasserot}). 
  If the group $G$ is understood, we drop it from the notation.
\end{Def}
\subsection{Properties of the Thom sheaf}
The following properties of the Thom sheaf are reformulations of 
well-known facts about Thom classes in
cohomology and $K$-theory. We deduce the elliptic case, whenever
the groups involved have been defined.

\medskip
\noindent
{\bf Naturality:}
Let $\map fXY$ be a $G$-equivariant map, and let $\xi$ be a
complex $G$-vectorbundle over $Y$. Then we have a commuting diagram of
sheaves over $X_{A_G}$ 
$$%
  \xymatrix{
    &\mathcal O_{X_{A_G}}\ar[1,-1]_{f_{A_G}^*\eta^\xi}\ar[1,1]^{\eta^{f^*\xi}}&\\
    \bL_G^{f^*\xi}\ar@{=}[0,2]^{\sim\phantom{xx}}&&   f_{A_G}^* \bL_G^\xi.
  }
$$
\begin{Pf}{}
  The map
  \begin{eqnarray*}
    \widetilde\mF_G\(f^\xi\)\negmedspace : \widetilde\mF_G\(Y^{\xi}\)
    &    \longrightarrow &
    \widetilde\mF_G\(X^{f^*\xi}\)
  \end{eqnarray*}
  is a map of $\mF_G(Y)$-module sheaves. Hence it corresponds to a map
  \begin{eqnarray*}
    \(\bL^{\xi}\)\inv & \longrightarrow & f_{A_G*}\(\bL^{f^*\xi}\)\inv,
  \end{eqnarray*}
  whose adjoint
    \begin{eqnarray*}
    f_{A_G}^*\(\bL^{\xi}\)\inv & \longrightarrow &\(\bL^{f^*\xi}\)\inv    
  \end{eqnarray*}
  is an isomorphism. In the elliptic case, the last statement follows
  from \cite[2.6]{Grojnowski}
  and the $W$-equivariance of the Thom isomorphism in cohomology.
\end{Pf}{}

\noindent
{\bf Change of groups:}
  Let $\map\phi HG$ be a map of Lie groups, and let $\map\xi VX$ be a
  $G$-equivariant complex vector bundle. Then we have a commuting
  diagram 
  of sheaves over $X_{A_{H}}$
  $$%
    \xymatrix{
      &\mathcal O_{X_{A_H}}\ar[1,-1]_{\eta_H^\xi}\ar[1,1]^{A_\phi^*\eta^{\xi}_G}&\\
      \bL_H^{\xi}\ar@{=}[0,2]^{\sim\phantom{xx}}&&   X_{A_\phi}^* \bL_G^\xi.
    }
  $$

  \begin{Pf}{}
    Consider
    the map of locally free rank one $\mF_H(X)$ module sheaves
    \begin{eqnarray*}
      \mF_H(X)\tensor_{A^*_\phi\mF_G(X)}
      A_\phi^*\widetilde\mF_G\(X^\xi\)
      &\longrightarrow&
      \widetilde\mF_H\(X^\xi\).
    \end{eqnarray*}
    For $[a]\in A_H$, we
    have an equality of fixed point sets
    \begin{eqnarray*}
      X^{T(a)} & = & X^{T(\phi(a))}.
    \end{eqnarray*}
    Hence \cite[2.6]{Grojnowski} implies that the above map 
    is an isomorphism at the stalk $[a]$.
  \end{Pf}

\noindent
{\bf Induction:}
  Assume that the induction axiom holds for the inclusion of a
  normal subgroup $K\vartriangleleft G$ and that $X$ is a $G$-complex on
  which the action of $K$ is free. Then the change of groups isomorphism
  $X_{A_{G}}\cong(K\backslash X)_{A_K}$ is covered by an isomorphism
  of line bundles identifying $\bL_G^\xi$ with
  $\bL_{G/K}^{K\backslash \xi}$ and $\eta_G^\xi$ with
  $\eta^{K\backslash\xi}_{G/K}$. 

\medskip
\noindent
{\bf Multiplicativity:} 
  Given equivariant complex vector bundles
  $\xi$ over a $G$-space $X$ and $\zeta$ over an $H$-space $Y$, 
  the K\"unneth isomorphism 
  $(X\times Y)_{A_{G\times H}}\cong X_{A_G}\times Y_{A_H}$
  is covered by an isomorphism of line bundles identifying 
  $
    \bL_{G\times H}^{\xi\oplus\zeta}
  $ 
  with the external tensor product 
  $
   \bL_G^{\xi}\tensor\bL_H^\zeta
  $
  and $\eta_{G\times H}^{\xi\oplus\zeta}$ with
  $\eta_G^\xi\tensor\eta_H^\zeta$. 
  In the special case where $G=H$, we get a commuting square of sheaves
  over $A_G$
  $$
    \xymatrix{
    &\mathcal O_{X_{A_G}}\ar[1,-1]\ar[1,1]& \\
    \bL_G^{\xi\oplus\zeta}\ar@{=}[0,2]^{\sim\phantom{XXX}}
    & &
   \bL_G^\xi\otimes_{\mathcal O_{A_G}}\bL_G^\zeta.    
    }
  $$

\begin{Pf}{}
The induction and (external) multiplicativity properties follow
directly from the induction axiom and the K\"unneth property of
$\mF_G$ in their formulation for pairs. The internal K\"unneth for
$\bL_T^{\xi\oplus \zeta}$ follows from the external K\"unneth and the
change of groups property.
\end{Pf}{}

\medskip
\noindent
{\bf Universal bundles:}
  The Thom sheaf of the universal complex line bundle
  $\xi_1$ over $\pt$ (compare Definition \ref{def:universal_bundle})
  is the line bundle
  \begin{eqnarray*}
    \bL_{U(1)}^{\xi_1} &\cong & \mL(0)
  \end{eqnarray*}
  of the divisor $(0)$ on
  $A_{U(1)}$. Recall that $A_{U(1)}$ equals $\CC$ or $\CC^\times$, in
  which case we
  replace $0$ by $1$, or $\CC/\langle\tau,1\rangle$.
  The Euler map of $\xi_1$ is the canonical inclusion
  \begin{eqnarray*}
    \mathcal O_{A_{U(1)}}&\longrightarrow & \mL(0).
  \end{eqnarray*}

  Consider the universal complex $n$-vector bundle
  $\xi_n$ over $\pt$. Let $T$ be the maximal torus of $U(n)$. Then
  \begin{eqnarray*}
    \bL_{T}^{\xi_n} &\cong & \bigotimes_{i=1}^n 
    p_i^*\medspace\mathcal L(0) 
  \end{eqnarray*}
  where 
  $$
    \longmap{p_i}{A_{U(1)}^n}{A_{U(1)}} 
  $$ 
  is the projection to the $i$th factor. 
 %
  This is the line bundle associated to the divisor 
  $$
    \sum_{t=i}^n \ker({p_i}),
  $$
  and the Euler map $\eta_T(\xi)$ is 
  its canonical inclusion of $\mathcal O_{A_{U(1)}^n}$ inside it.
  The $U(n)$-equivariant Thom sheaf and Euler map are obtained by
  taking the $S_n$-invariant parts of $\bL_T$ and $\eta_T$.
  We will write $\eta_n$ for the $n$th universal Euler map.

\medskip
The following result, which determines all Thom sheaves up to isomorphism,
follows immediately from the list of properties above. 
\begin{Thm}
  Let $\map\xi VX$ be a $G$-equivariant vector bundle, and
  let $c_\xi$ 
  be its Ginzburg-Kapranov-Vasserot characteristic class (c.f.\
  Definition \ref{def:characteristic_class}). 
  Then we have a commuting square
  $$
    \xymatrix{%
      \mathcal O_{X_{A_G}}\ar@{=}[r] \ar[d]_{\eta(\xi)} &
      {c_\xi}^*  \mathcal O_{A_{U(1)}^n/S_n} \ar[d]^{c_\xi^*(\eta_n)} \\ 
      \bL_G^{\xi}\ar@{=}[r]^{\sim\phantom{Xi}}&
      {c_\xi}^*  \bL^{\xi_{n}}.
   }
  $$  
\end{Thm}
\begin{Exa}\label{exa:e^l}
  Let $\lambda\neq 0$ be a weight of $T$, and let $j\negmedspace
  :K_\lambda\hookrightarrow T$ be the kernel of $e^{2\pi i\lambda}$. 
  Consider the
  $T$-equivariant line bundle 
  $\map{\xi_\lambda}{\CC_\lambda}\pt$.
  By Example \ref{exa:representations}, we have a short exact sequence
  $$
    A_{K_\lambda}\xrightarrow{\phantom{XX}A_j\phantom{XX}} A_T
    \xrightarrow{\phantom{XX}c_{\xi_{\lambda}}\phantom{XX}} A.
  $$
  Hence the Thom sheaf of $\xi_\lambda$ is 
  \begin{eqnarray*}
    \bL^{\xi_\lambda} & \cong & \mL(A_{K_\lambda}),
  \end{eqnarray*}
  the line bundle on $A_T$ associated to
  the divisor $A_j$.
\end{Exa}
\begin{Exa}\label{exa:representation_Thom_sheaf}
Let $\map\varrho GU(n)$ be a complex representation, viewed as a
$G$-equivariant complex vector bundle over the one point space. 
By Example \ref{exa:representations}, we have an isomorphism 
\begin{eqnarray*}
  \bL_G^{\varrho}&\cong& A_{\varrho}^*\medspace \bL_{U(n)}^{\xi_n}  
\end{eqnarray*}
of sheaves over $A_G$.  
\end{Exa}
\begin{Exa}
  \label{exa:trivial_representation}
  If $\map\varrho GU(1)$ is the one dimensional trivial
  represesntation, then $A_\varrho$ 
  factors through the inclusion of zero
  $\map{{i_0}}{A_1}A_{U(1)}$. Hence
  \begin{eqnarray*}
    \widetilde F_1(\SS^2) & \cong & i_0^* \medspace\mathcal I(0)
  \end{eqnarray*} is identified with the
  sheaf of invariant differentials on $A_{U(1)}$,
  \begin{eqnarray*}
    \omega &:=& \mathcal I(0)/\mathcal I(0)^2\at 0.
  \end{eqnarray*}
  Writing $\omega$ also for its pull-back to $A_G$, we obtain
  \begin{eqnarray*}
    \widetilde F_G(\SS^{2n}) & \cong & \omega^{\tensor n}.
  \end{eqnarray*} 
  This shows that the Periodicity Axiom (1.5.5) in
  \cite{Ginzburg:Kapranov:Vasserot} is, in part, redundant. We will
  come back to this in Section \ref{sec:coefficients}.
\end{Exa}
\begin{Exa}
  Let 
  $$
    \chi_\varrho = e^{2\pi i\lambda_1} + \dots + e^{2\pi i\lambda_n}
  $$
  be the character of $\varrho$, with $\lambda_k\in\Lambda\setminus\{0\}$ for 
  all $k$. 
  By Example \ref{exa:e^l}, we may identify
  the $G$-equivariant Thom sheaf 
  $\bL_G^{\varrho}$ over $A_T/W$ with the sheaf of $S_n$-invariant
  sections of
  $$
    \bigotimes_{i=0}^n \mathcal L(A_{K_{\lambda_i}}).
  $$
\end{Exa}

\begin{Cor}
  Let $\map\varrho TU(n)$ be a complex representation of $T$
  with character
  $$
    \chi_\varrho = e^{2\pi i\lambda_1} + \dots + e^{2\pi i\lambda n},
  $$
  $\lambda_i\neq 0$,
  and write $S_\varrho^{2n-1}$ for its unit sphere inside
  $\CC^n$. Then we have an isomorphism 
  \begin{eqnarray*}
    \(\SS_\varrho^{2n-1}\)_{A_T} & \cong &
    \bigcap_{i=1}^nA_{K_{\lambda_i}},
  \end{eqnarray*}
  where the right-hand side stands for the scheme theoretic
  intersection over $A_T$.
\end{Cor}
\begin{Pf}{}
  We have
  a cofiber sequence
  $$
    \(\SS^{2n-1}_\varrho\)_+\longrightarrow \SS^0 \longrightarrow
    \SS^{\varrho},  
  $$
  whose second map is the zero section $z$ in the definition of
  $\eta_T^{\varrho}$. Applying $\widetilde F_T$, we obtain the short exact
  sequence 
  $$
    \bigotimes_{i=1}^n\mathcal I(A_{K_{\lambda_i}})\longrightarrow
    \mathcal O_{A_{T}} 
    \longrightarrow 
    \mF_{T}\(\SS^{2n-1}_\varrho\),
  $$
  where the first map is the canonical inclusion.
\end{Pf}
\begin{Exa}\label{exa:infinitesimal}
  Consider the representation 
  \begin{eqnarray*}
    k\xi_1&:=&\xi_1\oplus\dots\oplus\xi_1
  \end{eqnarray*}
  of $U(1)$. Then 
  \begin{eqnarray*}
    \(\SS^{2k-1}_{k\xi_1}\)_{A_{U(1)}} & = & (0)^{[k]}
  \end{eqnarray*}
  is the $k$th infinitesimal neighbourhood of $0$ inside $A_{U(1)}$.  
\end{Exa}
We are now in a position to give the promised proof of the Completion
Theorem from the formal properties of $X_{A_T}$.

\medskip
\begin{Pf}{{of Theorem \ref{thm:completion}}}
  We follow the outline in
  \cite[(1.7.2)]{Ginzburg:Kapranov:Vasserot}. Writing $T$ as the
  product of $r$ copies of $U(1)$, we may build $ET$ from the
  equivariant skeleta
  \begin{eqnarray*}
    ET^{(2k-1)} & := & \SS^{2k-1}_{k\xi_1}\times\dots\times
    \SS^{2k-1}_{k\xi_1}, 
  \end{eqnarray*}
  where the notation is as in the last example. The induction axiom gives a
  commuting diagram
  $$%
    \xymatrix{
    \(ET^{(2k-1)}\times_T X\)_{A_1}\ar@{=}[0,2]^\sim\ar[d] &&
    \(ET^{(2k-1)}\times X\)_{A_T} \ar[d]\\
    A_1&&A_T.\ar[0,-2]
    }
  $$
  A combination of the internal and external K\"unneth properties,
  together with the last example,
  yields an isomorphism of schemes over $A_T$ 
  \begin{eqnarray*}
    \(ET^{(2k-1)}\times X\)_{A_T} &\cong &
    \( \prod_{j=1}^r\(0\)^{[k]} \)\times_{A_T} X_{A_T}    
  \end{eqnarray*}
  Letting $k$ vary, the right-hand side becomes an ind-scheme over $A_1$,
  isomorphic to the formal completion $\(X_{A_T}\)\widehat{{}_0}$. 
\end{Pf}
\subsection{$RO(G)$-grading and periodicity}
\label{sec:coefficients}
We are now ready to define the full theory $\Ell_G^*(-)$, graded by
the set of orthogonal representations contained in an indexing
universe (see \cite[p.154]{Alaska}) and their formal differences. For
any such universe, there is a cofinal system of representations of the
form 
$$
  \varrho\negmedspace : G\longrightarrow U(n) \longrightarrow O(2n).
$$
Hence it suffices to define the groups 
\begin{eqnarray*}
  \widetilde{\Ell}_G^{\varrho-\sigma}(X) 
  &:= &\bL^{\varrho}\tensor_{\mathcal O_{\mM_G}} \widetilde{\Ell}^0_G(\mathbb
  S^\sigma\wedge X),
\end{eqnarray*}
where $\varrho$ is as above and $\sigma$ is in our universe.
The resulting theory satisfies the axioms of a sheaf-valued
$RO(G)$-graded cohomology theory:
$\Ell_G$ is a contravariant functor of $X$ and $\sigma$ and a covariant
functor of $\varrho$. 
Each $\Ell_G^{\varrho-\sigma}(-)$ is exact on cofibre sequences and
sends wedges to products.
There are suspension isomorphisms
$$
  \isomap{s^\varrho}{\widetilde\Ell_G^{*+\varrho}(\mathbb S^\varrho\wedge
    X)}{\widetilde\Ell^*_G(X)}, 
$$
natural in $X$ and the orthogonal representation $\varrho$, and
satisfying 
$$
  s^{\varrho\oplus\sigma} =  s^{\sigma} \circ s^{\varrho}.
$$ 
\begin{Rem}
  These axioms are immediate from the definitions. 
  Note that exactness is checked on stalks, and we cannot expect
  sections over an open $\Gamma(U,\Ell_G^*(-))$ to be exact. In other
  words, a sheaf-valued cohomology theory is {\em not} the same thing as a sheaf
  of cohomology theories.
\end{Rem}
The theory of Thom sheaves is extended to virtual
equivariant complex vector bundles, and we set
\begin{eqnarray*}
  \Ell_G^{*+\xi}(X) & := & \widetilde\Ell_G^*(X^{-\xi}). 
\end{eqnarray*}
Finally, we have
periodicity isomorphisms 
\begin{eqnarray*}
  \widetilde{\Ell}_G^*(\mathbb S^\varrho\wedge X)  & \cong &
  \bL^{-\varrho}\tensor_{\Ell_G}\widetilde{\Ell}^*_G(X) \\
  &\cong & \widetilde{\Ell}_G^{*-\varrho}\tensor_{\Ell^*_G}\widetilde{\Ell}_G^*(X) 
\end{eqnarray*}
for complex representations
$\varrho$. 
\section{Euler classes}
Since the elliptic Thom sheaves are in general non-trivial, the notion
of Euler class does not have an immediate generalization to elliptic
cohomology. The different authors make different choices on this
matter, see \cite[p.10]{Rosu:elliptic}, \cite[5.3]{Ando:sigma} and
\cite[(2.6)]{Ginzburg:Kapranov:Vasserot}.  
\subsection{Thom isomorphisms}
Let $\xi_1$ be the universal complex line bundle of 
Definition \ref{def:universal_bundle}, and recall that its Thom sheaf
is
the invertible sheaf
\begin{eqnarray*}
  \bL^{\xi_1}_{U(1)}  & \cong & \mL(0)
\end{eqnarray*}
over $A_{U(1)}$.
For the additive or multiplicative group this divisor 
is principal:
\begin{eqnarray*}
  (0) & = & div(x) \quad\quad\text{on $\CC$ and} \\  
  (1) & = & div(1-z) \quad\quad\text{on $\CC^\times$}.  
\end{eqnarray*}
The {\em universal Euler classes} in cohomology and $K$-theory are the
functions 
$$%
  e(\xi_1)  =
  \begin{cases}
    x & \text{on $\CC$ and}  \\
    (1-z) & \text{on $\CC^\times$,}  
  \end{cases}
$$
and the 
{\em universal Thom isomorphisms} are
\begin{eqnarray*}
\vartheta\! :  \mathcal O_{A_{U(1)}} & \stackrel\cong\longrightarrow & \mL(0) \\
        f            & \longmapsto                   & \frac{f}{e(\xi_1)}
\end{eqnarray*}
(replace $(0)$ by $(1)$ for $K$-theory).
As a consequence, all Thom sheaves in these theories are trivialized,
and the theories possess Chern classes for complex vector bundles.

We will be particularly
interested in the case where the complex vector bundle $\map\xi VX$
comes equipped with a spin 
structure on the underlying real bundle. In this case, 
the Ginzburg-Kapranov-Vasserot characteristic class factors as
$$%
  \xymatrix{
    &A_{U^2(n)}\ar[1,1]^v& \\
    X_{A_G} \ar[0,2]_{c_\xi}\ar@{-->}[-1,1]^{\widetilde c_\xi} && A_{U(n)},
  }
$$
where $A_{U^2(n)}$ is the pull-back in the cartesian square
$$%
  \xymatrix{
    A_{U^2(n)}\ar[r]^u\ar[d]_v & A_{U(1)}\ar[d]^{\cdot 2}\\
    A_{U(n)}\ar[r]_{A_{\det}} & A_{U(1)}.
  }
$$
\begin{Exa}
  In the multiplicative case, a point in $A_{U^2(n)}$ consists of 
  $$
    (z_1,\dots,z_n)\in(\CC^\times)^n/S_n
  $$
  together with a choice of square root 
  $$
    (z_1\cdots z_n)^{\frac12}.
  $$  
\end{Exa}
\begin{Def}
  We write
  \begin{eqnarray*}
    \bL^{\xi_n}_{U^2(n)} & := & v^* \bL^{\xi_n}_{U(n)}.
  \end{eqnarray*}
  for the universal Thom sheaf for such $n$-dimensional complex bundles with
  spin structure.
\end{Def}
In $K$-theory $\bL^{\xi_n}_{U^2(n)}$ is the target of the {\em
  Atiyah-Bott-Shapiro Thom 
isomorphism}
\begin{eqnarray*}
\vartheta_{ABS}\!:
  \mathcal O_{A_{U^2(n)}} & \stackrel\cong\longrightarrow &
  \bL^{\xi_n}_{U^2(n)} \\ 
\notag        f            & \longmapsto                   & 
\frac{f}{e'(\xi_n)},
\end{eqnarray*}
with
\begin{eqnarray*}
  e'(\xi_n) & = &
  \(z_1^\frac12-z_1^{-\frac12}\)
  \cdots\(z_n^\frac12-z_n^{-\frac12}\).
\end{eqnarray*}
%
%
\subsection{Theta functions and elliptic Euler classes}
On the elliptic curve $E=\CC/2\pi i \langle
\tau,1\rangle$ the divisor $(0)$
is no longer principal. In this case, the Thom isomorphisms 
above are replaced by the theta function formalism:
let $\varrho$ be a complex representation of $G$
with character
$$e^{2\pi i\lambda_1}+\dots+e^{2\pi i\lambda_n},$$
$\lambda_k\in\Lambda\setminus\{0\}$.
Then the first Pontrjagin class
\begin{eqnarray*}
  p_1(\varrho)&=&\sum_{k=1}^n \lambda_k\tensor\lambda_k
\end{eqnarray*}
is an integer-valued positive definite symmetric bilinear form on
$\Lambda\!^\vee$. If $\varrho$ admits a spin structure then we have 
\begin{equation}\label{eq:2Z}
  p_1(\varrho) (x,x) \in 2\ZZ
\end{equation}
for all $x\in\Lambda\!^\vee$.   
\begin{Def}
  Let $I$ be a positive definite symmetric bilinear form on
  $\Lambda\!^\vee$, and assume that $I$ satisfies \eqref{eq:2Z}. Then  
  the {\em Looijenga line bundle} associated to $I$ is the invertible
  sheaf $\mL_I$ on
  on $\mM^h_T$ with sections
  \begin{eqnarray*}
    \mL_I(U) &=&
    \{
    f\in \Gamma\mathcal O^h_{y\inv U}   \mid
    f(q^xz) = q^{-\frac12I(x,x)}z^{-I^\sharp(x)}f(z) 
    \}. 
  \end{eqnarray*}
  Here $x\in\Lambda\!^\vee$, and $q^x$ stands for the image of $\tau x$ under  
  $\map\exp{\mathfrak t_\CC}T_\CC$, while $y$ is the quotient map 
  $T_\CC\to \mM_T$. So, if $T=U(1)$ then $q^x = e^{2\pi i\tau x}$.
  We write 
  \begin{eqnarray*}
    \Theta_\varrho & = & \Gamma \mL_{p_1(\varrho)}
  \end{eqnarray*}
  for the global sections of $\mL_{p_1(\varrho)}$ and refer to
  elements of $\Theta_\varrho$ as {\em Looijenga theta functions (of level $p_1(\varrho)$)}.
\end{Def}
\begin{Def}
  Let $\map\varrho GU(n)$ be as above. We are still assuming that
  we have a spin structure on $\varrho$.
  The {\em ``elliptic Euler class''} of $\varrho$ is the function on
  $T_\CC$ defined by
  \begin{eqnarray*}
    e_{ell}(\varrho) &:=& (-1)^n\prod_{i=1}^n\sigma\(q,z^{\lambda_i}\),
  \end{eqnarray*}
  where 
  \begin{eqnarray*}
    \sigma(q,z) & = &
    (z^{\frac12}-z^{-\frac12})\prod_{n\geq1}\frac{(1-q^nz)(1-q^nz\inv)}{\(1-q^{n}\)^2} 
  \end{eqnarray*}
  is the Weierstrass sigma function, and for $z=\exp(x)\in T_\CC$ and
  $\lambda\in\Lambda$ we are using the notation
  \begin{eqnarray*}
    z^{\lambda}\medspace&:=&\medspace e^{2\pi i\lambda(x)}.
  \end{eqnarray*}
\end{Def}
In elliptic cohomology, the role of the Thom isomorphism is replaced
by the isomorphism of 
line bundles over $\mM_G$
\begin{eqnarray*}
\vartheta\!:
  \mL_{p_1(\varrho)}^W & \stackrel\cong\longrightarrow &
  \bL^{\varrho}_{G} \\ 
\notag        f            & \longmapsto                   & 
\frac{f}{e_{ell}(\varrho)}.
\end{eqnarray*}
Presumably, these notions generalize to yield a theta
function description for the elliptic Thom
sheaf of any equivariant $U^2(n)$-bundle over a nice enough base (for
instance, an equivariantly formal space). We do not pursue this here. 
\begin{Rem}
  Ando has considered equivariant elliptic cohomology with twisted
  coefficients, where the twist comes from an element 
  $\beta\in H^4(BG;\ZZ)$. He does so by extending Looijenga's
  definition of $\mL_\beta$ to
  $\beta$ not satisfying \eqref{eq:2Z}. 
  Let $G$ be connected, and let $\varrho$ be an even dimensional
  orthogonal representation of $G$. 
  Then a similar argument to the one above yields an isomorphism between
  $\Ell^\varrho_G(\pt)$ and $\mL_{p_1(\varrho)}^W$. 
  So, the
  $RO(G)$-graded coefficients are contained in Ando's picture.
\end{Rem}

\subsection{Push-forwards}
\label{sec:push-forwards}
Let $X$ and $Y$ be compact, closed smooth manifolds, and let $\map fXY$ be a
complex oriented map in the sense of \cite{Quillen}. That means that
we have a factorization
$$
  \xymatrix{%
  &&E\ar[d]^\xi\\
  X\phantom{i}\ar@{>->}[-1,2]^i \ar[0,2]_f&&Y,
  }
$$ 
where $\xi$ is a complex vector bundle and the normal bundle $\nu$ of
$i$ is equipped with a complex structure. 
\begin{Def}
  For such a complex oriented map $f$, one defines the {\em relative
    Thom sheaf} as
  \begin{eqnarray*}
    \bL(f) &= & f_{A_G *}\bL^{-\nu}\tensor \bL^{\xi}.
  \end{eqnarray*}
  This is a sheaf over $Y_{A_G}$.  
  The {\em push-forward} along $f$ is the map
  \begin{eqnarray*}
    f_!\negmedspace:  \bL(f)
    & \longrightarrow& \mathcal O_{Y_{A_G}}
  \end{eqnarray*}
  of sheaves over $Y_{A_G}$ that is adjoint to the map
  \begin{eqnarray*}
    f_{A_G*}\bL^{-\nu} & \longrightarrow&
    \bL^{-\xi}
  \end{eqnarray*}
   induced by the Pontrjagin-Thom collapse.
\end{Def}
The following lemma is immediate from the definitions.
\begin{Lem}[Localization Lemma]
  Let $X$ be as above with a smooth
  $T$-action. Let 
  $i\negmedspace : 
  X^T\hookrightarrow X$ be the inclusion of the fixed points and
  assume that we are given a $T$-equivariant complex structure on the
  normal bundle 
  $\nu$ of $i$.
  Then we have a commuting diagram
  $$
    \xymatrix{
      &\mF_T^*(X)\ar[1,-1]_{i^*} & \\
      \mF_T^*(X^T)&&
      \mF_T^{*-\nu}(X^T),\ar[0,-2]_{z^*}\ar[-1,-1]_{i_!} \\
    }
  $$
  where $z$ is the zero section of $(X^T)^\nu$.
\end{Lem}
Note that the trivial representation does not turn up as a summand
inside $\nu$.
\begin{Cor}\label{cor:localization}
  In the situation of the Localization Theorem \ref{thm:localization},
  assume that the 
  normal bundle is equipped with a $T$-equivariant complex
  structure. Let $\Delta(\nu)\sub A_T$ be the closed subset 
  $$
    \Delta(\nu) = \bigcup_{\CC_\lambda\sub\nu} \ker\(A_{e^\lambda}\).
  $$
  Then, restricted to $A_T\setminus \Delta(\nu)$,
  the map $i^*$ becomes an isomorphism with inverse $i_!\circ(z^*)\inv$.
\end{Cor}

\section{Character Formulas}
\label{sec:character_formulas}
\subsection{Induced representations}
Let $G$ be a compact connected Lie group with maximal torus $T$, and
let $B\sub G_\CC$ be a Borel subroup of its complexification. Such a
choice of $B$ is equivalent to a choice of positive roots of
$G$. It endows the flag variety
\begin{eqnarray*}
  G/T&\cong& G_\CC/B  
\end{eqnarray*}
with a complex structure such that the tangent space at the coset of
$1$ is the complex $T$-representation
\begin{eqnarray*}
  \mathfrak g/\mathfrak t& \cong_\CC & \bigoplus_{\alpha\in\mR_-}\CC_\alpha.
\end{eqnarray*}
Similarly, if $H\sub G$ is a connected subgroup containing $T$ and $P_H$
the parabolic subgroup corresponding to $H$, we have a complex structure on
the homogenous space
\begin{eqnarray*}
  G/H &\cong& G_\CC/P_H.
\end{eqnarray*}
\begin{Def}\label{def:ind}
  In this situation, 
  we define the map
  $$
    \longmap{\on{ind}}{R(H)}{R(G)}
  $$
  as the composite
  $$
    \xymatrix{%
      K_H\ar@{=}[r]^{\sim\phantom{xxx}} \ar@{=}[d]_\vartheta & 
      K_G(G/H) \ar@{=}[d]^\vartheta \\
      K_H^{\mathfrak g/\mathfrak h}\ar@{=}[r]^{\sim\phantom{xxx}} &
      K_G^\tau(G/H)\ar[0,2]^{\pi_!}&&
      K_G.
    }
  $$
  Here, and I apologize for this notation, $\mathfrak h$ is the Lie
  algebra of $H$, {\em not} a Cartan subalgebra. Further 
  \begin{eqnarray*}
    \tau & \cong & G\times_H\mathfrak g/\mathfrak h
  \end{eqnarray*}
  is the tangent bundle
  of $G/H$, and $\pi$ is the unique map from $G/H$ to the one point
  space. 
  The push-forward
  $\pi_!$ is as defined in Section \ref{sec:push-forwards}.\footnote{A
    more common definition of the push-forward 
    $\pi_!$ in $K$-theory or cohomology is the composite of our
    $\pi_!$ with $\vartheta$. The reason for our convention is that it
    generalizes to elliptic cohomology in a canonical way.}
\end{Def}
The Atiyah-Singer index theorem identifies our definition of
$\on{ind}$ with the definition of induction found in the
representation-theory literature: 
\begin{Thm}[{c.f.\ \cite{Atiyah:Singer:ieoI} or \cite[5.4]{Hirzebruch:Berger:Jung}}]
  Let $\map\varrho HGL(V)$ be a complex representation. Then 
  \begin{eqnarray*}
    \on{ind}([\varrho]) & = & \sum (-1)^i H^i(G/H,\mathcal O(G\times_HV))    
  \end{eqnarray*}
  is the induced representation of $\varrho$. Here $\mathcal
  O(G\times_HV)$ is the sheaf of holomorphic sections of $G\times_HV$.
\end{Thm}
\subsection{The Weyl character formula}
We will now compute the character of these induced representations.
As in Theorem \ref{thm:flag}, we
let $W_H$ and $W_G$
be the respective Weyl groups,
and we let 
$$
  \longmap i{F}{G/H}
$$
be the inclusion of the $T$-fixed points $F:=(G/H)^T$. Recall that $F$
can be identified with the set $W_G/W_H$, and that we have
\begin{eqnarray*}
  i^*\tau & \cong_T & \coprod_{[w]\in F}(\mathfrak g/\mathfrak h)^w 
\end{eqnarray*}
(conjugation by $w$ on the right-hand side).
We have a commuting diagram
\begin{equation}
  \label{eq:weyl_diagram}
  \xymatrix{%
    K_H\ar@{=}[r]^{\sim\phantom{xxxx}} & 
    K_G(G/H)\ar@{=}[r]^\vartheta \ar[2,0]_{res}&
    K_G^\tau(G/H)\ar[0,2]^{\pi_!}\ar[2,0]_{res} & &
    K_G\ar[2,0]^{char}\\
    \\
    &
    K_T(G/H)\ar@{=}[r]^\vartheta \ar@{>->}[2,0]_{i^*}
    &
    K_T^\tau(G/H)\ar[0,2]^{\pi_!}\ar@{>->}[2,0]_{i^*}&&
    K_T\\
    \\
    &
    {\bigoplus\limits_{[w]\in F}K_T}\ar@{=}[r]^{\vartheta\phantom{xix}}&
    {\bigoplus\limits_{[w]\in F}K_T^{(\mathfrak g/\mathfrak h)^{w}}} &&
    {\bigoplus\limits_{[w]\in F}K_T}\ar@{>->}[0,-2]^{\phantom{xx}z^*}
    \ar[-2,-2]_{i_!}\ar[-2,0]_{\sum\limits_{[w]\in F}(-)_w},
  }  
\end{equation}
where $$\longmap zFF^{i^*\tau}$$ is the zero section.

The top row of \eqref{eq:weyl_diagram} is the map $\on{ind}$ of
definition \ref{def:ind}.
The composite of the vertical arrows on the left sends an
$H$-representation $\varrho$ to 
$
  \(\chi_{\varrho}^w\)_{[w]\in F}
$
(the character of $\varrho$ and its conjugates under $W_G$).
The composite at the bottom is multiplication by the Euler class
of $i^*\tau$. On the $[w]$th summand, this is
\begin{eqnarray*}
  e(i^*\tau)_{[w]}  & = & \prod_{\alpha\in\mR}\(1-e^{w(\alpha)}\).
\end{eqnarray*}
Here
\begin{eqnarray*}
  \mR &:= & \mR^G_-\setminus \mR^H_-
\end{eqnarray*}
consists of the negative roots of $G$ that are not roots of $H$.
Using the Localization Lemma (see Corollary \ref{cor:localization}),
we can deduce Weyl's character formula:
\begin{Thm}[Weyl]
  Let $\varrho$ be a representation of $H$. Then the character of its 
  induced representation equals
  \begin{eqnarray}
\label{eq:Weyl}    
\chi_{ind(\varrho)} & = & \sum_{[w]\in F}\medspace 
   \frac{\chi_\varrho^w}{\prod\limits_{\alpha\in\mR}\(1-e^{2\pi iw(\alpha)}\)}.
  \end{eqnarray}
\end{Thm}
To be precise, the Localization Lemma implies the equality
\eqref{eq:Weyl} in the localized ring 
$$
  R(T)[e(\mathfrak g/\mathfrak h)\inv].
$$
Since $R(T)$ maps injectively into this localization, and
$\chi_{ind(\varrho)}$ is an element of $R(T)$, it makes sense
to interpret \eqref{eq:Weyl} as a formula in $R(T)$.

Replacing $K$-theory by cohomology, we obtain a formula for the composite
$$
  \(Bj\)^*\circ \(Bk\)_!,
$$
where $j$ and $k$ are the respective inclusions of $T$ and $H$ in $G$.
Namely, it sends a regular function $f$ on $A_H$ to
\begin{eqnarray*}
  Bj^*(Bk_!(f)) & = &   \frac{\sum\limits_{[w]\in
      f}\det(w)f^w}{\prod\limits_{\alpha\in\mR}\alpha_\CC}. 
\end{eqnarray*}
\subsection{The Kac character formula}
We now turn our attention to the elliptic case, 
making the additional assumtion that the partial flag variety $G/H$
carries a $U^2$-structure. 
For simplicity of notation, we write
$$
  Ell_G^*(X) := \Gamma \Ell_G^*(X)^h
$$
for the analytic global sections, noting that the statement holds on
the level of sheaves with all sheaves pushed forward to $\mM_G$.
Let $\varrho$ be a $G$-representation.
The diagram \eqref{eq:weyl_diagram} is replaced by
\begin{equation*}
  \xymatrix{%
    \Theta_{\varrho+\mathfrak g/\mathfrak h}^{W_H}
    \ar@{=}[0,1]^{\sim\phantom{xxx}} &
    Ell_G^{\varrho+\tau}(G/H)\ar[0,2]^{\pi_!}\ar[2,0]_{res} & &
    Ell_G^\varrho\ar[2,0]^{char}
    \ar@{=}^\vartheta[r]&\Theta_{\rho}^{W_G}\ar[2,0]\\
    \\
(10)
&    Ell_T^{\varrho+\tau}(G/H)\ar[0,2]^{\pi_!}\ar@{>->}[2,0]_{i^*}&&
    Ell_T^\varrho\ar@{=}^\vartheta[r]&\Theta_{\rho}\\
    \\
    {\bigoplus\limits_{[w]\in F}\Theta_{\varrho+(\mathfrak g/\mathfrak
    h)^w}}\ar@{=}[0,1]^{\vartheta}
    &
    {\bigoplus\limits_{[w]\in F}Ell_T^{\varrho+(\mathfrak g/\mathfrak
        h)^{w}}} && 
    {\bigoplus\limits_{[w]\in F} Ell_T}^\varrho\ar@{>->}[0,-2]^{\phantom{xx}z^*}
    \ar[-2,-2]_{i_!}\ar[-2,0]_{\sum\limits_{[w]\in F}(-)_w}.   
  }  
\end{equation*}

\medskip
As before, we will write ${ind}$ for the composite of the arrows in
the top row. 
The Thom sheaf
$\bL^{\mathfrak g/\mathfrak h}_T$ is the line bundle over $\mM_T$
associated to the divisor 
$$
  \Delta = \sum_{\alpha\in\mathcal R} (A_{K_\alpha}),
$$
on $\mM_T$. Here $K_\alpha = \ker(e^{2\pi i\alpha})$.
%
The $T$-equivariant Euler class of $\mathfrak g/\mathfrak h$ is the
theta-function
\begin{eqnarray*}
  e_{ell}(\mathfrak g/\mathfrak h)&=&\pm
  \prod_{\alpha\in\mR} \(z^{\frac\alpha2}-z^{-\frac\alpha2}\)
  \prod_{n\geq 1}\frac{\(1-q^nz^\alpha\)\(1-q^nz^{-\alpha}\)}{(1-q^n)^2},
\end{eqnarray*}
where the sign equals $(-1)^{|\mR|}$. Set 
\begin{eqnarray*}
  \Phi\medspace=\medspace\Phi(q) &:=& \prod_{n\geq 1}\(1-q^n\)^2.
\end{eqnarray*}
\begin{Thm}\label{thm:Kac}
  Let $f$ be an element of $\Theta_{\varrho+\mathfrak g/\mathfrak
    h}^{W_H}$. Then we have 
  \begin{eqnarray*}
    ind(f) &=&\(-\Phi\)^{d}\frac
    {
    \sum\limits_{[w]\in F}\det(w)w(f)
    }
    {
    \prod\limits_{\alpha\in\mR} \(z^{\frac\alpha2}-z^{-\frac\alpha2}\)
    \prod\limits_{n\geq 1}\(1-q^nz^\alpha\)\(1-q^nz^{-\alpha}\)
    }.
  \end{eqnarray*}
\end{Thm}
Here $d$ is the complex dimension of the partial flag variety $G/H$.  
Consider now the special case where $G$ is simple and simply
connected, and $H=T$ is the maximal torus.
Then there is a smallest positive definite bilinear form $I_{Lo}$
satisfying \eqref{eq:2Z}. This is the bilinear form 
considered in \cite{Looijenga}, and we write $\mL_{Lo}$ for the
corresponding Looijenga line bundle.\footnote{In the notation of
  \cite{Looijenga} this is the line bundle $\mL\inv$.}
By \cite[(3.4)]{Looijenga}, we have 
\begin{eqnarray*}
  p_1(\mathfrak g/\mathfrak t) & = & g\cdot I_{Lo}
\end{eqnarray*}
and hence
\begin{eqnarray*}
  \mL_{p_1(\mathfrak g/\mathfrak t)}&=& \mL_{Lo}^g.
\end{eqnarray*}
Here $g$ is the dual Coxeter number. Assume that we have
$p_1(\varrho)=kI_{Lo}$ with $k\in\ZZ$. Then
the top row of (10) becomes a map
$$
  \longmap{ind}{\Theta_{k+g}}{\Theta_k^{W_G}},
$$
where $\Theta_k$ are the Looijenga theta functions of level $k$.
\begin{Def}[Looijenga basis]
  Let $k\in\mathbb N$,
  and let $\lambda\in \Lambda$. 
  The element $\theta_{k,\lambda}\in\Theta_k$
  is defined by 
  \begin{eqnarray*}
    \theta_{k,\lambda} &:=&
    \sum_{x\in\Lambda\!^\vee}q^{(k\phi+\lambda)(x)}e^{2\pi i(kI^\sharp(x)+\lambda)}.
  \end{eqnarray*}
  Here
  \begin{eqnarray*}
    \phi(x) &:=& \frac12I_{Lo}(x,x)
  \end{eqnarray*}
  and $\map{I^\sharp}{\Lambda\!^\vee}{\widehat T}$ is the adjoint of $I_{Lo}$.
\end{Def}
As $\lambda$ varies over a set of representatives for
$\Lambda/kI^\sharp(\Lambda\!^\vee)$, the $\theta_{k,\lambda}$ form a
basis for $\Theta_k$.  
\begin{Cor}[Kac character formula]
  In the situation of Theorem \ref{thm:Kac}, assume that $G$ is simple
  and simply connected and let $H=T$ be the maximal
  torus. Then we have
  \begin{eqnarray*}
    ind \(\theta_{k+g,\lambda}\) & = & \frac{\(-\Phi\)^{d}\cdot
    \sum\limits_{w\in W_G}\det(w)\cdot \theta_{k+g,w(\lambda)}
    }{\prod\limits_{\alpha\in\mR_-} \(e^{\pi i\alpha}-e^{-\pi i\alpha}\)
    \prod\limits_{n\geq 1}\(1-q^ne^{2\pi i\alpha}\)\(1-q^ne^{-2\pi i\alpha}\)}.  
  \end{eqnarray*}
\end{Cor}
Up to the factor 
$\pm\Phi(q)^{d+r}$, 
which is constant in $z$, this agrees with the Kac
character formula 
for the positive energy representation of the loop
group $\mL G$ of level $k$ and heighest weight
$$
  \lambda+\frac12\sum_{\alpha\in\mR_+}\alpha.
$$ 
For a presentation of the Kac character
  formula in this form see \cite[(14.3.4)]{Pressley:Segal} or \cite[11.4]{Ando:powerops}.
\bibliographystyle{alpha}
\bibliography{weyl.bib}
\end{document}